\documentclass[reqno,12pt,centertags]{amsart}
\usepackage{geometry}
\usepackage{amsmath,amsthm,amscd,amssymb,latexsym,upref,stmaryrd}
\usepackage{graphicx}
\usepackage{hyperref}
\newcommand*{\mailto}[1]{\href{mailto:#1}{\nolinkurl{#1}}}
\usepackage[numbers,sort&compress]{natbib}

\newcommand{\N}{\mathbb{N}}
\newcommand{\R}{\mathbb{R}}
\newcommand{\C}{\mathbb{C}}

\newcommand{\beq}{\begin{equation}}
\newcommand{\eeq}{\end{equation}}
\newcommand{\ba}{\begin{align}}
\newcommand{\ea}{\end{align}}

\newcommand{\dd}{{\rm d}}
\newcommand{\cC}{\mathcal C}

\newcommand{\cA}{\mathcal A}

\renewcommand{\ln}{\text{\rm ln}}
\newcommand{\ip}[2]{\left\langle #1,#2\right\rangle}

 \allowdisplaybreaks
\numberwithin{equation}{section}

\newtheorem{theorem}{Theorem}[section]
\newtheorem{lemma}[theorem]{Lemma}

\newtheorem{proposition}[theorem]{Proposition}

\theoremstyle{definition}

\newtheorem{remark}{Remark}[section]

\newtheorem{mainthm}{Theorem}

\newtheorem{romanthm}{Inverse Problem}

\begin{document}


\title[Half-inverse problems]
{Half-inverse problems for perturbed Bessel operators with a known left-hand potential}

\author[X.~J.~Xu]{Xin-Jian Xu}
\address{Department of Mathematics, School of Mathematics and Statistics, Nanjing University of
Science and Technology, Nanjing, 210094, Jiangsu, China}
\email{\mailto{xxjxuxinjian@163.com}}

\author[C.~F.~Yang]{Chuan-Fu Yang}
\address{Department of Mathematics, School of Mathematics and Statistics, Nanjing University of
Science and Technology, Nanjing, 210094, Jiangsu, China}
\email{\mailto{chuanfuyang@njust.edu.cn}}

\author[N. Bondarenko]{Natalia Bondarenko}
\address{Moscow Center of Fundamental and Applied Mathematics, Lomonosov Moscow State University, Moscow, Russia}
\email{\mailto{ bondarenkonp@sgu.ru}}

\author[R.~P.~Xia]{ Runpeng Xia}
\address{Department of Mathematics, Tandon School of Engineering, New York University, Brooklyn, NY 11201, USA}
\email{\mailto{rx2354@nyu.edu}}

\date{\today}

\begin{abstract}
{For the regular Sturm-Liouville operator $-\frac{{\rm d}^2}{{\rm d}x^2}+q$ acting in $L^2(0,1)$, the Hochstadt-Lieberman theorem states that one full spectrum together with the potential on either the left half-interval $q|_{(0,1/2)}$ or the right half-interval $q|_{(1/2,1)}$ uniquely determines the full  potential $q$. The inverse-square singular term $\tfrac{l(l+1)}{x^2}$ in the perturbed Bessel operator $-\frac{{\rm d}^2}{{\rm d}x^2}+\frac{l(l+1)}{x^2}+q$, acting in $L^2(0,1)$, breaks the left-right symmetry. It is known that, when the right half of the potential $q|_{(1/2,1)}$ is given, the uniqueness of the potential can be established even if a finite number of eigenvalues (depending on $l$) are omitted. In this paper we study the opposite setting--where the left half potential $q|_{(0,1/2)}$ is known--and show that one full spectrum together with finitely many norming constants uniquely determine the full potential $q$. Moreover, we specify the minimum number of norming constants required for the uniqueness. Finally, for positive integer $l$, $L^2$ potentials, and the Dirichlet boundary condition at $x = 1$, we establish the local bi-Lipschitz stability and the local solvability from the full spectrum and finitely many norming constants, assuming the potential on the left half-interval is fixed.}
\end{abstract}

\maketitle

\section{Introduction}

We consider the differential expression
\begin{equation}\label{1}
 \mathcal L_q:=-\frac{{\rm d}^2}{{\rm d}x^2} +\frac{l(l+1)}{x^2}+q(x),\quad x\in(0,1]. 
\end{equation}
Here $q\in L^1(0,1)$ is a real-valued potential and  $l\in[-\frac{1}{2},+\infty)$ is the angular-momentum quantum number. For $l\geq\frac{1}{2}$, the differential expression $\mathcal L_q$ is in the limit-point case at the singular endpoint $x=0$. For $-\frac{1}{2}\leq l<\frac{1}{2}$, however, it is in the limit-circle case and a usual boundary condition is imposed:
\begin{equation}\label{1.2}
\lim\limits_{x\rightarrow 0}x^{-(l+1)}f(x)<\infty.
\end{equation}
At the regular endpoint $x=1$, we impose
\begin{equation}\label{7}
 f'(1)+\beta f(1)=0,\quad \beta\in\mathbb{R}\cup\{\infty\},
\end{equation}
where $\beta=\infty$ means the Dirichlet condition $f(1)=0$. We denote by $H_{q,\beta}$ the self-adjoint operator in $L^2(0,1)$ generated by $\mathcal L_q$ subject to the boundary conditions \eqref{1.2} (when $-\frac{1}{2}\leq l<\frac{1}{2}$) and \eqref{7}. The boundary conditions and the self-adjointness of the operator $H_{q,\beta}$ can be found, e.g., in \cite{KST2010}. The operator $H_{q,\beta}$ is bounded below and has a simple discrete spectrum, which we denote by 
$$\sigma(q,\beta)=\{\lambda_n(q)|n\in\N\}.$$
Unless it would cause any ambiguity, we also write the  element $\lambda_n(q)$ as $\lambda_n$.

Let $\phi(x,\lambda)$ denote the solution of the equation $\mathcal L_q f=\lambda f$ normalized by
\begin{equation}\label{6}
 \lim_{x\rightarrow0}x^{-l-1}\phi(x,\lambda,q)=\frac{\sqrt\pi}{2^{l+1}\Gamma(l+\frac32)},
\end{equation}
where $\Gamma(\cdot)$ is the Gamma function. For $-\frac12\leq l<\frac12$, the relation \eqref{6} specifies the self-adjoint boundary condition \eqref{1.2} of $\mathcal L_q$ in the limit-circle case; for $l\geq\frac12$, it corresponds to a square-integrable solution. For $\lambda_n\in\sigma(q,\beta)$, the so-called norming constant is defined by 
\begin{equation}\label{11}
\tau_n(q):=\int_0^1\phi^2(x,\lambda_n, q)\dd x.
\end{equation}

For the perturbed Bessel operator and the related recent research, we refer to \cite{DF2023,GNZ2024,Shishkina2020} and the references therein. Inverse spectral problems for the operator $H_{q,\beta}$ have attracted considerable attention. For Borg-Levinson-type uniqueness theorems of the operator $H_{q,\beta}$, see, e.g., \cite{AHM2007, Carlson1997, Eckhart2014, FY2005, KST2010}. When the potential is known on a certain subinterval, uniqueness for the inverse problem solution can be found, e.g., in \cite{Eckhart2013, Guliyev2023, LSY2019, LX2026, XYYZ2023, DLL}. The problem of recovering the potential from eigenvalues associated with different angular-momentum quantum numbers $l$ has been studied, e.g., in \cite{GGHN, RS2001}. Local stability analysis for inverse problems can be found, e.g., in \cite{Serier2007}. Some reconstruction methods can be found, e.g., in \cite{AGM2008, AHM2008, HS2010}. The inverse problem theory for higher-order differential operators with Bessel-type singularities was developed in \cite{Yurko1993, Yurko1995} and for singular differential systems, in \cite{Ignatiev2023, Ignatiev2026}.

For a regular Sturm-Liouville operator on a finite interval, the corresponding regular half-inverse problem is classical. The Hochstadt-Lieberman theorem states that a full spectrum together with the potential on one half of the interval uniquely determines the full potential \cite{HL1978}; the known part of the potential may be given on either the left or the right half-interval. A variety of uniqueness results for partially given potential and spectral data were subsequently obtained by del Rio, Gesztesy, Simon, Horv\'ath, and others: the length of the unknown potential can be traded off against the density and type of the available spectral data, and optimal hypotheses can be expressed through closed exponential systems; see, e.g., \cite{dGS1997, GS2000, Horvath2001, Horvath2005, WX2009, WX2012, HS2016, Hat2021}. Generalizations of the Hochstadt-Lieberman theorem to quantum graphs are reviewed in \cite{Bondarenko2023}.

Unlike the regular Sturm-Liouville operator, the perturbed Bessel operator requires different spectral data for the uniqueness of the half-inverse problem solution when the potential is known on the left or right half of the interval, owing to the presence of the inverse-square singular term $l(l+1)/x^2$. To state the results, for $x\in\mathbb{R}$ we use the notations   
\begin{equation}\label{gauss}
\lceil x \rceil:={\rm min}\{n\in\mathbb{Z}|n\geq x\},\quad \lfloor x\rfloor:={\rm max}\{n\in\mathbb{Z}|n\leq x\}.
\end{equation}
And we use $\# A$ to represent the number of elements in the corresponding set $A$, counted with multiplicity when applicable. When the potential is known on the right half-interval $\left(\tfrac{1}{2},1\right)$,  a type of Hochstadt-Lieberman theorem (\cite[Theorem 4.3]{XYB2023}) can be stated as follows.  
\begin{mainthm}[\cite{XYB2023}]\label{xyb}
Let $l\in\N\cup\{0\}$, $p\in(1,+\infty)$, $\beta\in\mathbb{R}\cup\{\infty\}$, $q_1,q_2\in L^p(0,1)$, and $q_1=q_2$ a.e. on $\left(\tfrac{1}{2}, 1\right)$. Denote
\begin{equation}\label{nl}
n_l(\beta):=
\begin{cases}
\left\lfloor\frac{l}{2}\right\rfloor, \;\;\;\;\; \beta=\infty, \\
\left\lfloor\frac{l+1}{2}\right\rfloor,\;\; \beta\in\mathbb{R}.
\end{cases}
\end{equation}
If there exist two sets $A_1$ and $A_2$, such that $\# A_1=\# A_2=n_l(\beta)$ and
$$\sigma(q_1,\beta)\setminus A_1=\sigma(q_2,\beta)\setminus A_2,$$
then $q_1=q_2$ a.e. on $(0, 1)$.
\end{mainthm}
Theorem \ref{xyb} implies that, when the potential is known on the right half-interval, even an incomplete spectrum--missing finitely many eigenvalues (depending on $l$)--can still guarantee the uniqueness. Extensions to real $l\in [-\tfrac{1}{2}, +\infty)$, under additional smoothness assumptions near the midpoint, are considered in \cite[Corollary 5.4]{LX2026}. Some recent studies on the inverse problem of this type of right-part known potential $q|_{(a,1)}$ ($a\in(0,1)$) can be also found in \cite{DLL}.

The present paper treats the opposite situation: the potential is known on $(0,a)$ for some $a\in(0,1)$. The spectral data used in the inverse problem consist of partial information on eigenvalues and norming constants. To be precise, the inverse problem is formulated as follows: 
\begin{romanthm}\label{ip1}
Let $S\subset\Lambda\subset\sigma(q,\beta)$ and suppose that the potential $q(x)$ on $(0, a)$ is known a priori. Given the spectral data 
\begin{equation}\label{specdata}
\Lambda \;\;\text{and}\;\; \{\tau_n(q)|\lambda_n\in S\}, 
\end{equation}
find $q(x)$ on $(a, 1)$ and $\beta$.
\end{romanthm}

We apply the results of Horváth and Sáfár \cite{HS2016} to obtain closed-system criteria for the uniqueness of the solution to Inverse Problem \ref{ip1}. Theorem \ref{5.4} gives a necessary and sufficient condition in the Dirichlet case and a sufficient condition in the Robin case. In order to present the result, let $p\in[1,+\infty]$, $\frac{1}{p}+\frac{1}{p'}=1$, and let $I$ be a finite interval. 
A family $\mathcal E$ of functions in $L^{p'}(I)$ is called \emph{\textbf{closed}} in $L^p(I)$ if 
$$\text{for every}\; h\in L^p(I),\quad\int_I h(x)e(x)\,\dd x=0\quad(\forall e\in\mathcal E)\quad\Longrightarrow\quad h=0\; \text{a.e. on}\; I. $$
Define $b:=1-a$ and 
\begin{equation}\label{2}
 c_\lambda(y):=\cos(2\sqrt\lambda\,y),\quad
 d_\lambda(y):=-\frac{\partial}{\partial \lambda}c_\lambda(y)
\end{equation}
for $\lambda\in\C$ and $y\in (0,b)$; explicitly,
$$ d_\lambda(y)=\frac{y\sin(2\sqrt\lambda\,y)}{\sqrt\lambda}\quad(\lambda\neq0),\quad d_0(y)=2y^2.$$
For the sets $S\subset\Lambda\subset\sigma(q,\beta)$, let
\begin{equation}\label{3}
 \cC_b(\Lambda,S):=\{c_{\lambda_n}|\lambda_n\in\Lambda\}\cup\{d_{\lambda_n}|\lambda_n\in S\}.
\end{equation}

\begin{theorem}\label{5.4}
For fixed $a\in(0,1)$, $l\in[-\tfrac{1}{2},+\infty)$, and $p\in[1,+\infty]$, assume that $q_1,q_2\in L^1(0,1)$ and  
$$q_1=q_2\quad\text{a.e. on }(0,a),\quad q_2-q_1\in L^p(a,1).$$ 
Suppose that $\Lambda\subset\sigma(q_1,\beta_1)\cap\sigma(q_2,\beta_2)$ and $S\subset\Lambda$ satisfies 
$$\tau_n(q_1)=\tau_n(q_2), \quad \lambda_n\in S.$$

(1) If $\beta_1=\beta_2=\infty$, let $\mu\in\R$ with $\mu^2\notin\Lambda$, then $q_1=q_2$ a.e. on $(0,1)$
if and only if 
$\cC_b(\Lambda,S)\cup\{\cos(2\mu y)\}$ is closed in $L^p(0,b)$. 

(2) If $\beta_1, \beta_2\in\mathbb{R}$ and $\cC_b(\Lambda,S)$ is closed in $L^p(0,b)$, then $q_1=q_2$ a.e. on $(0,1)$ and $\beta_1=\beta_2$.
\end{theorem}

The only if part of Theorem \ref{5.4} (1) means that  if $ \cC_b(\Lambda,S)\cup\{\cos(2\mu y)\}$ is not closed in $L^p(0,b)$, then there exist $q_1$ and $q_2$ satisfying all the conditions in the theorem with $q_1\neq q_2$ in $L^1(0,1)$.

When the potential is known on more than the left half-interval, uniqueness from a single spectrum follows from the general Hochstadt-Lieberman-type theorem of Eckhardt and Teschl \cite[Theorems 5.3 and 5.4]{Eckhart2013}. In the setting of the present paper, we state it as the following theorem. 

\begin{mainthm}[\cite{Eckhart2013}]\label{ET2013}
Suppose that $l\in[-\frac12, +\infty)$, $\beta_1,\beta_2\in\mathbb{R}\cup\{\infty\}$, and $q_1,q_2\in L^1(0,1)$. If there exists $\varepsilon\in\left(0,\tfrac12\right)$ such that 
$$q_1=q_2\quad\text{a.e. on }\left(0,\tfrac12+\varepsilon\right),\quad \text{and}\quad \sigma(q_1,\beta_1)=\sigma(q_2,\beta_2),$$
then $q_1=q_2$ a.e. on $(0,1)$ and $\beta_1=\beta_2$.
\end{mainthm}

We now consider the critical case in which the potential is known on exactly the interval $\left(0,\tfrac12\right)$. The following Theorem \ref{7.2} gives explicit sufficient numbers of additional norming constants for the uniqueness.

\begin{theorem}\label{7.2}
Suppose that $l\in[-\frac12, +\infty)$, $p\in[1,+\infty]$, and $q_1,q_2\in L^1(0,1)$, with $\left(q_1-q_2\right)|_{(\tfrac{1}{2},1)}\in L^p(\tfrac{1}{2},1)$. Assume
$$q_1=q_2\quad\text{a.e. on }\left(0,\tfrac12\right),\quad \text{and}\quad \sigma(q_1,\beta_1)=\sigma(q_2,\beta_2).$$
Let $S\subset\sigma(q_1,\beta_1)$ such that 
$$\tau_n(q_1)=\tau_n(q_2), \quad \lambda_n\in S.$$
Then $q_1=q_2$ a.e. on $(0,1)$ and $\beta_1=\beta_2$ under either of the following sufficient conditions: for $p\in[1,+\infty)$,  
\begin{equation}\label{36}
\# S \geq M_D^{(p)}(l) :=\max\!\left\{0, \left\lceil\frac{l-2+1/p}{2}\right\rceil\right\},\quad \beta_1=\beta_2=\infty, 
\end{equation}
\begin{equation}\label{35}
\# S\geq M_R^{(p)}(l):=\max\!\left\{0, \left\lceil\frac{l-1+1/p}{2}\right\rceil\right\},\quad \beta_1, \beta_2\in\mathbb{R},
\end{equation}
and for $p=+\infty$, 
\begin{equation}\label{36bb}
\# S \geq M_D^{(\infty)}(l) :=\max\left\{0,\left\lfloor\frac {l}{2}\right\rfloor\right\},\quad \beta_1=\beta_2=\infty, 
\end{equation}
\begin{equation}\label{35bb}
\# S\geq M_R^{(\infty)}(l):=\max\left\{0,\left\lfloor\frac{l+1}{2}\right\rfloor\right\},\quad \beta_1, \beta_2\in\mathbb{R}.
\end{equation}
\end{theorem}

In the following Theorems \ref{7.3} and \ref{7.0}, we show that the numbers $M_D^{(p)}(l)$ and $M_R^{(p)}(l)$ in Theorem \ref{7.2} are uniformly minimal for the uniqueness. For this purpose, for every model potential $q\in L^p(0,1)$, we construct a distinct potential that agrees with $q$ on the left half-interval and shares its full spectrum and any prescribed $M_D^{(p)}(l)-1$ or $M_R^{(p)}(l)-1$ norming constants, respectively. The Robin parameter is allowed to vary in Theorem \ref{7.0}. The proofs combine explicit constructions based on Poisson integral representations of Bessel functions \eqref{i38}, the Paley-Wiener theorem, and an adaptation of Horv\'ath's nonlinear mapping method (see \cite{Horvath2005}). 

\begin{theorem}\label{7.3}
Let $l\in[-\frac12, +\infty)$, $p\in[1,+\infty]$, $q\in L^p(0,1)$, $M_D^{(p)}(l)\geq1$, and $S\subset\sigma(q,\infty)$ with $\# S= M_D^{(p)}(l)-1$.  
Then there exists $\tilde q\in L^p(0,1)$ such that
$$\tilde q=q\quad\text{a.e. on }(0,\tfrac12),\quad \tilde q\neq q, \quad \sigma(\tilde q,\infty)=\sigma(q,\infty),$$
and
$$\tau_n(\tilde q)=\tau_n(q),\quad \lambda_n\in S.$$
\end{theorem}

\begin{theorem}\label{7.0}
Let $l\in[-\frac12, +\infty)$, $p\in[1,+\infty]$, $\beta\in\R$, $q\in L^p(0,1)$, $M_R^{(p)}(l)\geq1$, and $S\subset\sigma(q,\beta)$ with $\# S=M_R^{(p)}(l)-1$. Then there exist $\hat{q}\in L^p(0,1)$ and $\hat{\beta}\in\R$ such that 
\begin{equation}\label{r42}
\hat{q}=q\quad\text{a.e. on }(0,\tfrac12),\quad \hat{q}\neq q, \quad \sigma(\hat{q},\hat{\beta})=\sigma(q,\beta),
\end{equation}
and
\begin{equation}\label{r43}
\tau_n(\hat{q})=\tau_n(q),\quad \lambda_n\in S.
\end{equation}
\end{theorem}

By comparing Theorem \ref{xyb} with Theorems \ref{7.2}--\ref{7.0}, we observe that, for the Hochstadt-Lieberman-type problems in which the potential is known on the left or right half-interval, the spectral data required for the uniqueness are significantly different. More interesting, if $l\in\N\cup\{0\}$, then $M_D^{(p)}(l)=\left\lfloor\frac{l}{2}\right\rfloor$ and $M_R^{(p)}(l)=\left\lfloor\frac{l+1}{2}\right\rfloor$, which are $n_l(\beta)$ for $\beta=\infty$ and $\beta\in\R$ defined in \eqref{nl}, respectively.


Savchuk and Shkalikov \cite{SS2010, SS2013} establish general results on the uniform stability and solvability of inverse Sturm-Liouville problems in Sobolev scales. Stability of half-inverse problems for Sturm-Liouville operators on a finite interval, particularly the Hochstadt-Lieberman problem, has been studied in, e.g., \cite{Bondarenko2025, HK2010, YB2020}; solvability conditions and constructive methods have been obtained in, e.g., \cite{MP2010, Sakh2001}.

Basing on the uniqueness result in Theorem \ref{7.2}, we investigate local stability and solvability of the corresponding half-inverse problem. For positive integer $l$ and $L^2$ potentials, we consider the Dirichlet problem with a fixed potential on the interval $\left(0,\tfrac{1}{2}\right)$. The spectral data consist of the full spectrum and $M_D^{(2)}(l)$ selected norming constants. We show that the norming data generate derivative functions which, together with the spectral cosine functions and the constant function, form a Riesz basis. This yields bounded invertibility of the linearized spectral map (see Lemma \ref{isom}). Based on this map, we obtain the local bi-Lipschitz stability and the local solvability of the inverse problem.

To formulate these results, we fix a model potential $q^*\in L^2(0,1)$ and put
\begin{equation}\label{41}
 \mathcal X_{q^*}:=\{q\in L^2(0,1) | q=q^* \;\text{a.e. on }(0,\tfrac{1}{2})\}.
\end{equation}
By identifying $\mathcal X_{q^*}$ with $L^2(\tfrac{1}{2},1)$, we consider the map 
$$\mathcal S_J: \mathcal X_{q^*}\left(\cong L^2\left(\tfrac{1}{2},1\right)\right)\longrightarrow \mathcal Y:=\R\oplus\ell^2(\N)\oplus\R^M$$ 
defined by 
\begin{equation}\label{44}
 \mathcal S_J(q)=\big(m(q-q^*),\{\lambda_n(q)-\lambda_n(q^*)-m(q-q^*)\}_{n\geq1},\{\tau_j(q)-\tau_j(q^*)\}_{j\in J}\big),
\end{equation}
where $J\subset\N$ with $M:=\# J=M^{(2)}_D(l)$, and $m(\cdot)$ is defined by 
\begin{equation}\label{43}
m(f):=\int_{1/2}^1f(x)\,\dd x,\; f\in L^2(\tfrac{1}{2},1).
\end{equation}
The map $\mathcal S_J$ is well defined since for $l\in\mathbb{N}$ and $q\in L^2(0,1)$, the eigenvalues $\lambda_n(q)\in\sigma(q,\infty),\; n\in\mathbb{N},$ satisfy
\begin{equation}\label{45}
 \lambda_n(q)=\pi^2\left(n+\frac l2\right)^2 +\int_0^1q(x)\,\dd x-l(l+1)+s_n(q),\quad \{s_n(q)\}_{n\geq1}\in\ell^2.
\end{equation}
It is well known that the sequence $\{s_n(q)\}_{n\geq1}$ in \eqref{45} is an $\ell^2(\mathbb{N})$-valued real-analytic function of $q$, and each fixed $\tau_j(q)$ is real analytic as well (see \cite{Serier2007}). Thus, the map $\mathcal S_J:\mathcal X_{q^*}\rightarrow\mathcal Y$ is real analytic. Theorem \ref{7.2} implies that the map $\mathcal S_J$ is injective. Moreover, since $q^*$ is fixed, from the asymptotic estimate \eqref{45}, the value $\mathcal S_J(q)$ in $\mathcal{Y}$ is equivalent to the spectral data $\{\lambda_n(q)\}_{n\geq1}$ and $\{\tau_j(q)\}_{j\in J}$.

For $q_1,q_2\in\mathcal X_{q^*}$, we define 
\begin{equation}\label{49}
 \begin{aligned}
 d_J\left(\mathcal S_J(q_1), \mathcal S_J(q_2)\right)^2&:= \|\mathcal S_J(q_1)-\mathcal S_J(q_2)\|^2_{\mathcal Y}  \\
 &=|m(q_2-q_1)|^2+\sum_{n=1}^{\infty}\left|\lambda_n(q_2)-\lambda_n(q_1)-m(q_2-q_1)\right|^2\\
 &+\sum_{j\in J} \left|\tau_j(q_2)-\tau_j(q_1)\right|^2.
 \end{aligned}
\end{equation}


The following theorem establishes the local bi-Lipschitz stability for the half-inverse problem. 
\begin{theorem}\label{8.4}
Let $l\in\N$, $\beta=\infty$, $q^*\in L^2(0,1)$, $J\subset\N$ with $\# J=M^{(2)}_D(l)$. Then there exist constants $C_*, C_1, C_2>0$, depending on $l$, $q^*$, and $J$, such that, whenever
$$q_i\in\mathcal X_{q^*},\quad \|q_i-q^*\|_{L^2(1/2,1)}<C_*,\quad i=1,2,$$
one has
\begin{equation}\label{98}
C_1\|q_2-q_1\|_{L^2(1/2,1)}\leq d_J\left(\mathcal S_J(q_1),\mathcal S_J(q_2)\right) \leq C_2\|q_2-q_1\|_{L^2(1/2,1)}.
\end{equation}
\end{theorem}

The following theorem establishes the local solvability of the half-inverse problem.
\begin{theorem}\label{8.6}
Let $l\in\N$, $\beta=\infty$, $q^*\in L^2(0,1)$, $J\subset\N$ with $\# J=M^{(2)}_D(l)$. Then there exist constants $\delta>0$ and $C>0$, depending only on  $l$, $q^*$, and $J$, such that, for every $\mathfrak{s}\in\mathcal Y$ satisfying 
$$\|\mathfrak{s}\|_{\mathcal Y}<\delta,$$
there exists a unique potential $q\in  \mathcal X_{q^*}$ such that 
$$\mathcal S_J(q)=\mathfrak{s}.$$
Moreover, 
\begin{equation}\label{th162}
\|q-q^*\|_{L^2(1/2,1)}\leq C \|\mathfrak{s}\|_{\mathcal{Y}}.
\end{equation}
\end{theorem}

The paper is organized as follows. Section \ref{sec2} provides the proofs of Theorems \ref{5.4} and \ref{7.2}. In Section \ref{sec3}, we construct an entire function $F(\lambda)$, which is used in the proof of Theorem \ref{7.0} and in the local stability and solvability analysis. 
Section \ref{sec4} provides the proofs of Theorems \ref{7.3} and \ref{7.0}. Section \ref{local} provides the proofs of Theorems \ref{8.4} and \ref{8.6}. Appendix \ref{appe} contains two auxiliary lemmas. 
\vskip 1cm

\section{Proofs of Theorems \ref{5.4} and \ref{7.2}}\label{sec2}

In this section, by using the results of \cite[Theorems 1.11-1.13]{HS2016}, we give the proofs of Theorems \ref{5.4} and \ref{7.2}.

\begin{proof}[Proof of Theorem \ref{5.4}.]
 For $i=1,2$, denote 
$$V_i(x):=\frac{l(l+1)}{x^2}+q_i(x),\quad x\in (a,1).$$
Since $a>0$, we have $V_i\in L^1(a,1)$ and $V_2-V_1=q_2-q_1\in L^p(a,1)$. 
The normalization condition \eqref{6} and $q_1=q_2$ almost everywhere on $(0,a)$ imply
\begin{equation}\label{restc}
\phi(x,\lambda,q_1)=\phi(x,\lambda,q_2)=:\phi(x,\lambda),\quad 0<x\leq a.
\end{equation}

It can be seen that every $\lambda_n\in\Lambda$ is also an eigenvalue of the corresponding Sturm-Liouville problem 
\begin{equation}\label{slp}
\begin{cases}
    & -f''(x)+V_i(x)f(x)=\lambda f(x),\quad x\in[a,1], \\
     & \phi(a,\lambda_n)f'(a)-\phi'(a,\lambda_n)f(a)=0,\quad f'(1)+\beta_i f(1)=0
\end{cases}
\end{equation} 
for both $i=1,2$. 

For $\lambda_n\in S$, the norming constants of the Sturm-Liouville problems \eqref{slp} corresponding to $\lambda_n$ are given by
$$T_{i,n}:=\frac{\tau_n(q_i)-\displaystyle\int_0^a\phi(x,\lambda_n)^2\,dx}{\phi(a,\lambda_n)^2+\phi'(a,\lambda_n)^2}$$
for $i=1,2$ respectively. Then $\tau_n(q_1)=\tau_n(q_2)$ gives $T_{1,n}=T_{2,n}$ for $\lambda_n\in S$.

The proof is complete by applying Theorems 1.11 and 1.12 in \cite{HS2016} to the Sturm-Liouville problems \eqref{slp}. 
\end{proof}

\begin{remark}\label{5.2} 
Theorem \ref{5.4} (2) provides only a sufficient criterion for uniqueness. The necessity direction requires additional hypotheses, which are not imposed here. Further details on the hypotheses for the regular Sturm-Liouville operator setting can be found in \cite{Horvath2005}. 
\end{remark}

To prove Theorem \ref{7.2}, we need the following asymptotic estimates for the eigenvalues. For $q\in L^1(0,1)$ and $l\in[-\tfrac{1}{2},+\infty)$, the eigenvalues $\lambda_{n}\in\sigma(q,\beta)$ with $n\in\mathbb{N}$ have the following asymptotic behavior as $n\rightarrow \infty$ (see \cite[Theorem 2.5]{KST2010}):  
\begin{equation}\label{as}
\sqrt{\lambda_{n}}=
\begin{cases}
\left(n+\dfrac{l}{2}\right)\pi+r_n(\beta), & \beta=\infty,\\
 \left(n+\dfrac{l-1}{2}\right)\pi+r_n(\beta), & \beta\in\mathbb{R}, 
\end{cases}
\end{equation}
where
\begin{equation}\label{asr}
r_n(\beta)=
\begin{cases}
 O\left(\frac{1}{n}\right), & l>-\frac{1}{2}, \\
O\left(\frac{\ln n}{n}\right), & l=-\frac{1}{2}. 
\end{cases}
\end{equation}

Before giving the proof of Theorem \ref{7.2} by applying \cite[Theorem 1.13]{HS2016}, we now establish a class of integrated counting asymptotics. For a discrete set $\Gamma\subset\mathbb{R}$, let
$$n_\Gamma(t):=\#\{\lambda\in\Gamma| |\lambda|\leq t\},\quad t\geq 0.$$

\begin{lemma}\label{7.1}
Suppose that $0\notin\sigma(q,\beta)$. Then for large $R$ one has  
\begin{equation}\label{33}
 \int_1^R\frac{2n_{\sigma(q,\beta)}(t^2)}t\,\dd t=
 \begin{cases}
 \frac{2R}{\pi}-(l+1)\ln R+O(1),&\beta=\infty,\\
 \frac{2R}{\pi}-l\ln R+O(1),&\beta\in\R.
 \end{cases}
\end{equation}
If $S\subset\sigma(q,\beta)$ consists of $M$ distinct fixed eigenvalues, then
\begin{equation}\label{34}
 \int_1^R\frac{2n_S(t^2)}t\,\dd t=2M\ln R+O(1).
\end{equation}
\end{lemma}

\begin{proof} 
If $l>-\frac{1}{2}$, by using \cite[Lemma 3.2]{HS2016}, the asymptotic formula \eqref{as} gives \eqref{33}. 

If $l=-\frac{1}{2}$, the asymptotic formula \eqref{as} gives
$$\sqrt{\lambda_{n}}=
\begin{cases}
\left(n-\frac{1}{4}\right)\pi+r_n(\beta), & \beta=\infty,\\
 \left(n-\frac{3}{4}\right)\pi+r_n(\beta), & \beta\in\mathbb{R}, 
\end{cases}$$
with $r_n(\beta)=O\left(\frac{\ln n}{n}\right)$. 
The argument in the proof of \cite[Lemma 3.2]{HS2016} remains applicable. More precisely, 
denote 
$$\sqrt{\lambda'_{n}}:=
\begin{cases}
\left(n-\frac{1}{4}\right)\pi, & \beta=\infty, \\
\left(n-\frac{3}{4}\right)\pi, & \beta\in\mathbb{R},
\end{cases}$$
and $\sigma_{\beta}:=\{\lambda'_n|n\in\mathbb{N}\}$. Then $\int_1^R\frac{2n_{\sigma_{\beta}}(t^2)}t\,\dd t$ satisfies the same asymptotic estimate as \eqref{33}. Moreover,
\begin{equation}\label{bb33}
\begin{split}
\left|\int_1^R\frac{2n_{\sigma(q,\beta)}(t^2)-2n_{\sigma_{\beta}}(t^2)}{t}\,\dd t \right|& \leq C_0+ 2\sum_{n=N}^{\infty}\left|\ln\frac{\sqrt{\lambda'_n}+r_n(\beta)}{\sqrt{\lambda'_n}}\right|  \\
     & \leq C_0+ C\sum_{n=1}^{\infty}  \frac{\ln n}{n^2} <\infty
\end{split}
\end{equation}
for some constant $C>0$, where $N\geq1$ is chosen so that $\lambda_n,\lambda_n'>1$ for all $n\ge N$, and the contributions of the finitely many terms, including any negative eigenvalues, are absorbed into $C_0$. Thus the perturbation \eqref{bb33} can be absorbed into the $O(1)$ remainder in \eqref{33}, and the conclusion still holds. 

Let 
$$ R_0:=\max\left\{1,\max_{\lambda\in S}\sqrt{|\lambda|}\right\}. $$
For $R\geq R_0$, one has $n_S(t^2)=M$ for $t\ge R_0$, and hence 
$$ \int_1^R\frac{2n_S(t^2)}{t}\,dt = \int_1^{R_0}\frac{2n_S(t^2)}{t}\,dt +2M\int_{R_0}^R\frac{dt}{t} =2M\ln R+O(1). $$
This proves \eqref{34}. 
\end{proof}

Now we give the proof of Theorem \ref{7.2}.

\begin{proof}[Proof of Theorem \ref{7.2}]
After a common spectral shift, if necessary, we may suppose that all eigenvalues are positive.  Such a shift changes both potentials by the
same constant, preserves their equality on the known interval, leaves the eigenfunctions and the norming constants unchanged, and does not affect uniqueness. 

We first consider the case $p\in[1,+\infty)$. Assume that $\#S=M$. Put 
$$I(R):=\int_1^R\frac{2n_{\sigma(q_1,\beta_1)}(t^2)+2n_S(t^2)}{t}\,dt.$$
By Lemma \ref{7.1}, for large $R$, one has 
\begin{equation}\label{ir1}
I(R)=
\begin{cases}
\frac{2R}{\pi}+(2M-l-1)\ln R+O(1),&\beta_1=\beta_2=\infty,\\
\frac{2R}{\pi}+(2M-l)\ln R+O(1),&\beta_1,\beta_2\in\mathbb{R}.
\end{cases}
\end{equation}
Applying \cite[Theorem 1.13]{HS2016},  one has the uniqueness when 
\begin{equation}\label{ir2}
\limsup_{R\to\infty}\left[I(R)-\frac{2R}{\pi}+\left(2+\frac1{p'}\right)\ln R\right]>-\infty,\quad \beta_1=\beta_2=\infty,
\end{equation}
or
\begin{equation}\label{ir3}
\limsup_{R\to\infty}\left[I(R)-\frac{2R}{\pi}+\frac1{p'}\ln R\right]>-\infty,\quad \beta_1,\beta_2\in\mathbb{R}.
\end{equation}
Here $1/p+1/p'=1$. 

The preceding asymptotics \eqref{ir1} show that the condition \eqref{ir2} or \eqref{ir3} holds whenever
$$2M-l+1+\frac1{p'}\ge0,\quad \text{for}\quad \beta_1=\beta_2=\infty,$$
or
$$2M-l+\frac1{p'}\ge0,\quad \text{for}\quad \beta_1,\beta_2\in\mathbb{R}.$$
Since $M$ is a nonnegative integer and $1/p+1/p'=1$, these are exactly the bounds \eqref{36} and \eqref{35}, respectively. Consequently, $q_1=q_2$ almost everywhere on $(0,1)$, and $\beta_1=\beta_2$.

It remains to consider the case $p=\infty$. For each fixed $l$, the definitions of the counting numbers in \eqref{36}-\eqref{35bb} imply that, for all sufficiently large finite $s$, 
$$M_D^{(s)}(l)=M_D^{(\infty)}(l),\quad M_R^{(s)}(l)=M_R^{(\infty)}(l).$$
Since
$$(q_1-q_2)|_{(1/2,1)}\in L^\infty\left(\tfrac{1}{2},1\right)\subset L^s\left(\tfrac{1}{2},1\right),$$
the conclusion follows by applying the case $p\in[1,+\infty)$, with $p$ replaced by $s$. This completes the proof.
\end{proof}

\vskip 1cm
\section{Common spectral data function $F(\lambda)$}\label{sec3}

This section introduces an auxiliary entire function $F$ encoding common eigenvalues and norming constants, and derives its cosine integral representations. These results are used in the proof of Theorem \ref{7.0} in Section \ref{sec4} and the local stability and solvability analysis in Section \ref{local}.

In addition to the solution $\phi(x,\lambda, q)$, we denote another solution of the equation $\mathcal L_q f=\lambda f$ by $\psi(x,\lambda, q)$, which satisfies \eqref{7} and the following initial conditions
\begin{equation}\label{8}
 \begin{aligned}
 &\psi(1,\lambda, q)=0,\quad
   \psi'(1,\lambda, q)=-1,
   &&\beta=\infty,\\
 &\psi(1,\lambda, q)=1,\quad
   \psi'(1,\lambda, q)=-\beta,
   &&\beta\in\R.
 \end{aligned}
\end{equation}
For every fixed $x>0$, the solutions $\phi(x,\lambda, q)$ and $\psi(x,\lambda, q)$ are entire functions of $\lambda$; see, e.g., \cite{KST2010}.  For a fixed $\beta$ under consideration, we can construct the so-called characteristic function
\begin{equation}\label{9}
 \Delta(\lambda,q):=W\!\left(\psi(\cdot,\lambda, q), \phi(\cdot, \lambda, q)\right), 
\end{equation}
where $W(\cdot,\cdot)$ is the Wronskian and it is defined by 
$$W(f,g):=fg'-f'g.$$ 
It is easy to see that 
\begin{equation}\label{deldef}
\Delta(\lambda,q)=
\begin{cases}
\phi(1,\lambda,q),&\beta=\infty,\\
\phi'(1,\lambda,q)+\beta\phi(1,\lambda,q), &\beta\in\R.
\end{cases}
\end{equation}
The function $\Delta(\lambda, q)$ is entire in $\lambda$, and its zeros correspond to the spectrum  
$$\sigma(q,\beta)=\{\lambda_n|n\in\mathbb{N}\}.$$

At an eigenvalue $\lambda_n\in\sigma(q,\beta)$, there is a nonzero constant $\kappa_n(q)$ such that
\begin{equation}\label{10}
 \psi(x,\lambda_n, q) =\kappa_n(q)\phi(x,\lambda_n, q).
\end{equation}
By considering the norming constant $\tau_n(q)$ defined in \eqref{11}, one has (see, e.g., (2.16) in \cite{XYYZ2023})
\begin{equation}\label{12}
 \dot\Delta(\lambda_n, q) =-\kappa_n(q)\tau_n(q),
\end{equation}
where the dot denotes differentiation with respect to $\lambda$, that is, $\dot\Delta(\lambda, q):=\frac{\partial}{\partial\lambda}\Delta(\lambda, q)$.

Assume that 
$$q_1=q_2\quad \text{a.e. on}\; (0,a),$$ 
and fix the right boundary parameters $\beta_1,\beta_2\in\mathbb{R}\cup\{\infty\}$. We compare the two operators 
$$H_{q_i,\beta_i},\quad i=1,2,$$ 
with the same angular-momentum quantum number $l\in[-\frac{1}{2},+\infty)$. 

When no ambiguity arises, we write
$$\phi_i(x,\lambda):=\phi(x,\lambda,q_i),\quad \psi_i(x,\lambda):=\psi(x,\lambda,q_i),\quad \Delta_i(\lambda):=\Delta(\lambda, q_i),\quad i=1,2.$$
Define
\begin{equation}\label{15}
 F(\lambda)
 :=W\!\left(\psi_2(\cdot, \lambda),\psi_1(\cdot, \lambda)\right)(a).
\end{equation}
This is an entire function of $\lambda$. 

\begin{lemma}\label{3.1}
Suppose $l\in[-\frac12, +\infty)$, $q_1,q_2\in L^1(0,1)$, $a\in (0,1)$, $q_1=q_2$ a.e. on $(0,a)$. If $\lambda_n\in\sigma(q_1,\beta_1)\cap\sigma(q_2,\beta_2)$ for $\beta_1, \beta_2\in\R$ or $\beta_1=\beta_2=\infty$, then 
\begin{equation}\label{16}
F(\lambda_n)=0, 
\end{equation}
and
\begin{equation}\label{17}
 \dot F(\lambda_n)=\kappa_n(q_1)\kappa_n(q_2)\bigl(\tau_n(q_1)-\tau_n(q_2)\bigr). 
\end{equation}
\end{lemma}

\begin{proof}
The normalization condition \eqref{6} and $q_1=q_2$ almost everywhere on $(0,a)$ imply
\begin{equation}\label{14}
\phi_1(x,\lambda)=\phi_2(x,\lambda)=:\phi(x,\lambda),\quad 0<x\leq a.
\end{equation}
If $\lambda_n\in\sigma(q_1,\beta_1)\cap\sigma(q_2,\beta_2)$, then \eqref{14} gives
$$\psi_i(x,\lambda_n)=\kappa_{n}(q_i)\phi(x,\lambda_n),\quad 0<x\leq a,\quad i=1,2.$$
Hence the Wronskian in \eqref{15} vanishes, and therefore $F(\lambda_n)=0$. From \eqref{10}, \eqref{14}, and \eqref{15}, by direct calculation we obtain 
\begin{equation}\label{c1}
\dot F(\lambda_n)=\kappa_n(q_1)W\left(\dot\psi_2(\cdot,\lambda_n),\phi(\cdot,\lambda_n)\right)(a)-\kappa_n(q_2)W\left(\dot\psi_1(\cdot,\lambda_n),\phi(\cdot,\lambda_n)\right)(a).
\end{equation} 
Note that for $\lambda=\lambda_n$, we have
\begin{equation}\label{c2}
\kappa_n(q_1)W\left(\psi_2(\cdot,\lambda_n),\dot{\phi}(\cdot,\lambda_n)\right)(a)-\kappa_n(q_2)W\left(\psi_1(\cdot,\lambda_n),\dot{\phi}(\cdot,\lambda_n)\right)(a)=0.
\end{equation} 
The two identities \eqref{c1} and \eqref{c2} give   
$$\dot F(\lambda_n)=\kappa_n(q_1)\dot\Delta_2(\lambda_n)-\kappa_n(q_2)\dot\Delta_1(\lambda_n).$$
By considering \eqref{12} one has \eqref{17}.
\end{proof}

\begin{remark}\label{3.2}
The converse statement needed later is also important.  Suppose that $\lambda_n\in\sigma(q_1,\beta_1)$ and $F(\lambda_n)=0$. Then the two vectors $\bigl(\psi_i(a,\lambda_n),\partial_x\psi_i(a,\lambda_n)\bigr)$, $i=1,2$, are  proportional.  Since the equations coincide on $(0,a)$, there exists a constant $c\neq0$, such that 
$$\psi_2(x,\lambda_n)=c\phi_1(x,\lambda_n),\quad x\in(0,a].$$
Hence $\lambda_n$ is an eigenvalue of $H_{q_2,\beta_2}$ with the \emph{same Bessel condition} \eqref{1.2} at $x=0$. That is,
$$\lambda_n\in\sigma(q_2,\beta_2).$$
\end{remark}

Remember that $b=1-a$. For $i=1,2$, set
\begin{equation}\label{18}
 y:=1-x,\quad \quad u_i(y,\lambda):=\psi_i(1-y,\lambda),
\end{equation}
and
\begin{equation}\label{19}
 Q_i(y):=\frac{l(l+1)}{(1-y)^2}+q_i(1-y), \quad 0<y<b.
\end{equation}
Then
$$-u_i''+Q_i u_i=\lambda u_i$$
are two regular Sturm-Liouville equations on $(0,b)$. Moreover, the solutions $u_i(y,\lambda)$ satisfy the initial conditions:
$$u_i(0,\lambda)=1,\quad u_i'(0,\lambda)=\beta_i,\quad \beta_i\in\mathbb{R},$$
or
$$u_i(0,\lambda)=0,\quad u_i'(0,\lambda)=1, \quad \beta_i=\infty.$$
We define 
\begin{equation}\label{20}
 h(y):=Q_2(y)-Q_1(y)=q_2(1-y)-q_1(1-y).
\end{equation}
The singular Bessel term cancels in the difference. In view of the above, all transformation-operator arguments below take place on a regular interval.

\begin{proposition}\label{4.1}
Let $a\in (0,1)$, $q_j\in L^1(0,1)$, $q_1=q_2$ a.e. on $(0,a)$, and $h\in L^p(0,b)$ defined by \eqref{20} with $p\in[1,+\infty]$.

(1) If $\beta_1=\beta_2=\infty$, then for every $\mu\in\mathbb{R}$, there is an operator 
$$A_{\mu}=A_{\mu}(Q_1,Q_2): L^p(0,b)\rightarrow L^p(0,b)$$ 
depending on $Q_1,Q_2$, such that 
\begin{equation}\label{22}
-2(z^2-\mu^2)F(z^2) =-\int_0^b h(y)\,\dd y  +\int_0^b\cos(2zy)(A_{\mu}h)(y)\,\dd y. 
\end{equation} 
In particular, $A_{\mu}$ is an isomorphism on $L^p(0,b)$.

(2) If $\beta_1,\beta_2\in\R$, then there is an operator 
$$B=B(Q_1,\beta_1,Q_2,\beta_2): L^p(0,b)\rightarrow L^p(0,b)$$ 
depending on $Q_1,\beta_1,Q_2,\beta_2$,  such that
\begin{equation}\label{z21}
 F(z^2)=\beta_2-\beta_1+\frac{1}{2}\int^b_0h(y)\,\dd y+ \int_0^b\cos(2zy)(Bh)(y)\,\dd y.
\end{equation}
In particular, $B$ is an isomorphism on $L^p(0,b)$.
\end{proposition}

\begin{proof} 
From \eqref{15} and \eqref{18}-\eqref{20}, using the Lagrange identity one has
\begin{equation}\label{fz}
F(z^2)=
\begin{cases}
  \int_{0}^{b}h(y)u_1(y,z^2)u_2(y,z^2)\,\dd y, & \beta_1=\beta_2=\infty,\\
  \beta_2-\beta_1+\int_{0}^{b}h(y)u_1(y,z^2)u_2(y,z^2)\,\dd y, & \beta_1,\beta_2\in\R.
\end{cases}
\end{equation}

(1) From \eqref{fz} with $\beta_1=\beta_2=\infty$, and the product representation of $u_1(y,z^2)u_2(y,z^2)$ in \cite[Lemma 5.2]{Horvath2005}, one can conclude \eqref{22}. Moreover, 
$$A_{\mu}=I+V_{\mu},$$
where $V_{\mu}$ is a Volterra operator on $L^p(0,b)$. Thus $A_{\mu}$ is an isomorphism on $L^p(0,b)$.

(2) 
Using \eqref{fz} with $\beta_1,\beta_2\in\R$, the product representation of $u_1(y,z^2)u_2(y,z^2)$ in \eqref{aa2}, and Fubini's theorem, one can obtain \eqref{z21}, where the representation of $B$ is given in \eqref{aa3}. The  representation \eqref{aa3} also shows that $B$ is an isomorphism on $L^p(0,b)$.
\end{proof}

\vskip 1cm
\section{Proofs of Theorems \ref{7.3} and \ref{7.0}}\label{sec4}

Next we will prove that $M_D^{(p)}(l)$ and $M_R^{(p)}(l)$ in Theorem \ref{7.2} are the minimum numbers of norming constants required for the uniqueness. That is, we prove Theorems \ref{7.3} and \ref{7.0}.

\begin{proof}[Proof of Theorem \ref{7.3}]
We first consider the case $p\in[1,+\infty)$. 
For $\lambda=z^2$, the solution $\phi(x,z^2, q)$ with $q=0$ is
\begin{equation}\label{ii38}
\phi(x,z^2, 0)=z^{-l-1}\sqrt{\frac{\pi zx}{2}}J_{l+\frac12}(zx).
\end{equation}
From Poisson's integral representation for the Bessel function (see (9.1.20) in \cite{AS1972}) 
\begin{equation}\label{i38}
J_{l+1/2}(z)=\frac{1}{2^{l-1}\sqrt{2\pi}\Gamma(l+1)}z^{l+\tfrac{1}{2}}\int_{0}^{1}(1-t^2)^l\cos(zt)\,\dd t,
\end{equation}
together with \eqref{deldef}, one has
\begin{equation}\label{38}
 \Delta(z^2,0)=\phi(1,z^2, 0) =\frac{1}{2^l\Gamma(l+1)}\int_0^1(1-t^2)^l\cos(zt)\,\dd t.
\end{equation}

For convenience, let $M_D:=M_D^{(p)}(l)\geq1$. 
Choose $\mu\neq0$ with $\mu^2\notin\sigma(q,\infty)$ and set
$$P(\lambda):=(\lambda-\mu^2) \prod_{\lambda_n\in S}(\lambda-\lambda_n),\quad
 w_l(t):=(1-t^2)^l,\quad 0<t<1.$$
The empty product is understood to be $1$. Since $\# S=M_D-1$, we have $\deg P=M_D$. The definition of $M_D$ gives $M_D<\frac{1}{2}\left(l+1/p\right)$, and hence $l>2M_D-1/p$. 
For $0\leq j\leq2M_D$, we have
$$w_l^{(j)}(t)=O\left((1-t)^{l-j}\right),\quad t\rightarrow 1^-.$$
Consequently, $w_l\in W^{2M_D,\;p}(0,1)$. Moreover, since 
$$l>2M_D-\frac{1}{p}\geq2M_D-1,$$
the endpoint values satisfy
$$w_l^{(j)}(1)=0, \quad 0\leq j\leq2M_D-1.$$
Every odd derivative of $w_l$ vanishes at $0$, because $w_l$ is even in a neighborhood of $0$. Thus
$$h_D:=P\!\left(-\frac{\dd^2}{\dd t^2}\right)w_l\in L^p(0,1),$$
and repeated integration by parts produces no boundary terms.
Formula \eqref{38} yields
\begin{equation}\label{39}
 \int_0^1h_D(t)\cos(zt)\,\dd t=2^l\Gamma(l+1)P(z^2)\Delta(z^2,0).
\end{equation}

To extend the representation \eqref{39} to the model potential $q$, write 
$$E(z):=P(z^2)\bigl(\Delta(z^2,q)-\Delta(z^2,0)\bigr).$$
Since $q\in L^1(0,1)$, the estimate in \cite[Lemma 2.2]{KST2010} together with \eqref{deldef} gives
$$\left|\Delta(z^2,q)-\Delta(z^2,0)\right|\leq\frac{C_q {\rm e}^{|\operatorname{Im}z|}}{(1+|z|)^{l+2}}$$
for some $C_q>0$, and therefore
$$|E(z)|\leq C_{q,S}(1+|z|)^{2M_D-l-2}{\rm e}^{|\operatorname{Im}z|}$$
for some $C_{q,S}>0$. Since 
$$2M_D-2-l<1/p-2\leq -1,$$
the even entire function $E$ has exponential type at most $1$ and belongs to $L^1(\mathbb R)\cap L^2(\mathbb R)$ on the real axis. By the  Paley-Wiener theorem, its inverse Fourier transform is supported in $[-1,1]$; its $L^1$ integrability also ensures that this inverse transform is bounded and continuous. Since $E$ is real and even on the real axis, there exists a real-valued function $h_1\in C[0,1]$ such that
\begin{equation}\label{xxez}
E(z)=\int_0^1h_1(t)\cos(zt)\,dt.
\end{equation}
Consequently,
$$\tilde h_D:=h_D+2^l\Gamma(l+1)h_1\in L^p(0,1)$$
satisfies
\begin{equation}\label{xx39}
\int_0^1\tilde h_D(t)\cos(zt)\,dt=2^l\Gamma(l+1)P(z^2)\Delta(z^2,q).
\end{equation}

The right-hand side in \eqref{xx39} is not identically zero, so $\tilde h_D\neq0$. With $t=2y$ and $\lambda=z^2$, the relation \eqref{xx39} becomes 
\begin{equation}\label{40}
\int_0^{1/2}\tilde h_D(2y)\cos(2\sqrt\lambda\,y)\,\dd y=2^{l-1}\Gamma(l+1)P(\lambda)\Delta(\lambda,q).
\end{equation}
This yields
\begin{equation}\label{b3}
\begin{cases}
\int_0^{1/2} \tilde h_D(2y)c_{\lambda_n}(y)\,\dd y=0, & \lambda_n\in\sigma(q,\infty), \\
  \int_0^{1/2} \tilde h_D(2y)d_{\lambda_n}(y) \,\dd y=0, &   \lambda_n\in S,\\
  \int_0^{1/2} \tilde h_D(2y)c_{\mu^2}(y)\,\dd y=0. & 
\end{cases}
\end{equation}
The identities \eqref{b3} imply $\cC_{1/2}(\sigma(q,\infty),S)\cup\{\cos (2\mu y)\}$ is not closed in $L^p(0,\tfrac{1}{2})$.  From Theorem \ref{5.4} (1) there exists $\tilde q\neq q$  such that
$$\tilde q=q\quad\text{a.e. on }(0,\tfrac12),$$
and every $\lambda_n\in\sigma(q,\infty)$ remains an eigenvalue corresponding to $\tilde q$, and the norming constants agree on $S$.  Thus
$$\sigma(q,\infty)\subset\sigma(\tilde q,\infty).$$
Then the asymptotic behavior of the eigenvalues \eqref{as} implies 
$$\sigma(q,\infty)=\sigma(\tilde q,\infty).$$

It remains to consider the case $p=+\infty$. Put $M_D:=M_D^{(\infty)}(l)\geq1$ in this case. Since $l\geq 2M_D$, we have 
$$w_l\in W^{2M_D,\infty}(0,1), \quad w_l^{(j)}(1)=0,\quad 0\le j\le2M_D-1.$$
Together with the vanishing odd derivatives at zero, these identities justify the same integrations by parts. Moreover, $2M_D-l-2\leq-2$, so the Paley-Wiener argument remains valid and yields $\tilde h_D\in L^\infty(0,1)$. Thus \eqref{39}-\eqref{b3} hold as before. Then Theorem \ref{5.4} (1) and the eigenvalue asymptotics \eqref{as} give the required potential $\tilde q\in L^\infty(0,1)$. This completes the proof.
\end{proof}

\begin{proof}[Proof of Theorem \ref{7.0}]
We first consider the case $p\in[1,+\infty)$. Set $M_R:=\# S=M^{(p)}_R(l)-1\geq0$. The definition of $M^{(p)}_R(l)$ gives 
\begin{equation}\label{r44}
 M_R= \left\lceil\frac{l-1+1/p}{2}\right\rceil-1<\frac{1}{2}\left(l-1+1/p\right).
\end{equation}
Consequently, $l>2M_R+1-1/p$.
Recall that \eqref{deldef} for $\beta\in\R$ gives 
$$\Delta(\lambda,q)=\phi'(1,\lambda,q)+\beta\phi(1,\lambda,q).$$ 
Equations \eqref{ii38} and \eqref{i38} imply 
 $$\phi(x,z^2,0)=\frac{1}{2^l\Gamma(l+1)}x^{l+1}\int_0^1(1-t^2)^l\cos(zxt)\,\dd t.$$
Differentiating with respect to $x$, putting $x=1$, and integrating by parts once, we obtain  
\begin{equation}\label{r45}
 \Delta(z^2,0)=\frac{1}{2^l\Gamma(l+1)}\int_0^1\omega_{l,\beta}(t)\cos(zt)\,\dd t,
\end{equation}
where
$$\omega_{l,\beta}(t):=(1-t^2)^{l-1}\bigl[(l+\beta)+(l-\beta)t^2\bigr].$$

Set  
$$P_S(\lambda):=\prod_{\lambda_n\in S}(\lambda-\lambda_n),\quad h_R(t):=P_S\left(-\frac{\dd^2}{\dd t^2}\right)\omega_{l,\beta}(t),$$
with the empty product understood to be $1$. Since $\# S=M_R$, we have $\deg P_S=M_R$. For $0\leq j\leq2M_R$, there holds 
$$\omega_{l,\beta}^{(j)}(t)=O\bigl((1-t)^{l-1-j}\bigr),\quad t\rightarrow1^-.$$
Since $l-1-2M_R>-1/p$, we have
$$\omega_{l,\beta}\in W^{2M_R,\;p}(0,1),\quad h_R\in L^p(0,1).$$
Moreover, for $0\leq j\leq 2M_R-1$, 
$$l-1-j\geq l-2M_R>0,$$
and hence 
$$\omega^{(j)}_{l,\beta}(1)=0.$$
Since $\omega_{l,\beta}$ is even in a neighborhood of zero, every odd derivative of $\omega_{l,\beta}$ vanishes at $0$. Thus repeated integration by parts in \eqref{r45} produces no boundary terms.  It yields
$$\int_0^1h_R(t)\cos(zt)\,\dd t=2^l\Gamma(l+1)P_S(z^2) \Delta(z^2,0).$$

To extend this representation to the model potential $q$, define 
$$E_{\beta}(z):=P_S(z^2)\bigl(\Delta(z^2,q)-\Delta(z^2,0)\bigr).$$
Since $q\in L^1(0,1)$, the estimates in \cite[Lemma 2.2]{KST2010} together with \eqref{deldef} imply 
$$\left|\Delta(z^2,q)-\Delta(z^2,0)\right|\leq\frac{C_{q,\beta}{\rm e}^{|\operatorname{Im}z|}}{(1+|z|)^{l+1}}$$
for some $C_{q,\beta}>0$. Therefore,
$$|E_{\beta}(z)|\leq C_{q,\beta,S}(1+|z|)^{2M_R-l-1}{\rm e}^{|\operatorname{Im}z|}$$
for some $C_{q,\beta,S}>0$. Since 
$$2M_R-l-1<\frac1p-2\leq -1,$$
the same argument used to obtain \eqref{xxez} yields a real-valued function $h_2\in C[0,1]$ such that
$$E_{\beta}(z)=\int_0^1h_2(t)\cos(zt)\,dt.$$
Consequently, the function
$$\tilde h_R :=h_R+2^l\Gamma(l+1)h_2\in L^p(0,1)$$
satisfies
\begin{equation}\label{inthr}
\int_0^1\tilde h_R(t)\cos(zt)\,dt=2^l\Gamma(l+1)P_S(z^2)\Delta(z^2,q).
\end{equation}
The right-hand side of \eqref{inthr} is not identically zero, so $\tilde h_R\neq0$.

For a nonzero $\gamma\in\R$, we define $g_{\gamma}(y):=\gamma \tilde h_R(2y)$ for $0<y<\tfrac12$.  The substitution $t=2y$ together with $\lambda=z^2$ then gives 
\begin{equation}\label{r46}
 \int_0^{1/2}g_{\gamma}(y)c_\lambda(y)\,\dd y=2^{l-1}\Gamma(l+1)\gamma P_S(\lambda)\Delta(\lambda,q).
\end{equation}
The right-hand side of \eqref{r46} vanishes at every $\lambda\in\sigma(q,\beta)$ and has a double zero at every $\lambda\in S$. Hence 
\begin{equation}\label{r3}
\begin{cases}
\int_0^{1/2} g_{\gamma}(y)c_{\lambda_n}(y)\,\dd y=0, & \lambda_n\in\sigma(q,\beta), \\
  \int_0^{1/2} g_{\gamma}(y)d_{\lambda_n}(y) \,\dd y=0, &   \lambda_n\in S.
\end{cases}
\end{equation}

In the construction leading to \eqref{15}, set $q_1(x)=q(x)$, $\beta_1=\beta$, and 
$$ q_2(x)=
\begin{cases}
q(x),&0<x<\tfrac{1}{2},\\[1mm]
q(x)+r(1-x),&\tfrac{1}{2}<x<1,
\end{cases}\quad \beta_2=\beta+\delta,$$
where $r\in L^p(0,\frac12)$ and $\delta\in\R$. Then, we denote $F$ in \eqref{15} by $F_{r,\delta}$. Let  
$$Q(y):=\frac{l(l+1)}{(1-y)^2}+q(1-y), \quad 0<y<\frac12.$$
By \eqref{z21} we may define the operator 
$$B(r,\delta):=B(Q,\beta, Q+r,\beta+\delta): L^p(0,\tfrac{1}{2})\rightarrow L^p(0,\tfrac{1}{2}),$$
which depends on $Q,\beta,Q+r,\beta+\delta$, such that 
\begin{equation}\label{r47}
 F_{r,\delta}(z^2)
 =\delta+\frac12\int_0^{1/2}r(y)\,\dd y
  +\int_0^{1/2}\cos(2zy)
   [B(r,\delta)r](y)\,\dd y.
\end{equation}

Let $\tilde r\in L^p(0,\frac12)$, $\tilde\delta\in\R$. From Lemma \ref{a.1}, on every bounded set of parameters, there exists a constant $C>0$, such that  
\begin{equation}\label{r48}
 \|B(r,\delta)-B(\tilde r,\tilde\delta)\|\leq C\bigl(\|r-\tilde r\|_{L^p(0,1/2)}+|\delta-\tilde\delta|\bigr). 
\end{equation}
Set
$$
 \delta(r):=-\frac12\int_0^{1/2}r(y)\,\dd y.$$
Hölder's inequality gives 
\begin{equation}\label{holder}
|\delta(r)-\delta(\tilde r)|\leq\frac{1}{2}\left(\frac{1}{2}\right)^{1-1/p}\|r-\tilde r\|_{L^p(0,1/2)}.
\end{equation}
Consequently, \eqref{r48} and \eqref{holder} show that the map from $r$ to $B(r,\delta(r))$ is locally Lipschitz in the operator norm. Moreover, $B(0,0)$ is an isomorphism by Proposition \ref{4.1}. Therefore, Lemma \ref{h2005} in Appendix \ref{appe} implies that, for every small enough $\gamma\ne0$,  there exists $r_\gamma\in L^p(0,\frac12)$ such that
\begin{equation}\label{r49}
 B(r_\gamma,\delta(r_\gamma))r_\gamma=g_\gamma.
\end{equation}
Also we get $r_\gamma\ne0$ since $g_\gamma\ne0$.

Define
$$\hat{q}(x)=
\begin{cases}
q(x),&0<x<\frac12,\\
q(x)+r_\gamma(1-x),&\frac12<x<1,
\end{cases}
\quad \hat{\beta}=\beta+\delta(r_\gamma).$$
Then \eqref{r47}, \eqref{r49}, and the definition of $\delta(r)$ give
\begin{equation}\label{r50}
F_{r_\gamma,\delta(r_\gamma)}(\lambda)
 =\int_0^{1/2}g_\gamma(y)c_\lambda(y)\,\dd y.
\end{equation}
By \eqref{r46}, we obtain 
$$F_{r_\gamma,\delta(r_\gamma)}(\lambda_n)=0\quad(\lambda_n\in\sigma(q,\beta)),
 \quad
 \dot F_{r_\gamma,\delta(r_\gamma)}(\lambda_n)=0\quad(\lambda_n\in S).$$
Lemma \ref{3.1} and Remark \ref{3.2} then imply   
\begin{equation}\label{r5100}
\sigma(q,\beta)\subset\sigma(\hat{q},\hat{\beta}), 
\end{equation}
and
\begin{equation}\label{r51}
 \tau_n(\hat{q})=\tau_n(q),\quad \lambda_n\in S.
\end{equation}
It follows from the eigenvalue asymptotics \eqref{as} and \eqref{r5100} that 
\begin{equation}\label{r52}
\sigma(\hat{q},\hat{\beta})=\sigma(q,\beta). 
\end{equation}
Equations \eqref{r51} and \eqref{r52} imply \eqref{r42} and \eqref{r43}. Moreover, $\hat{q}\neq q$ since $r_\gamma\neq0$, and clearly $\hat{q}\in L^p(0,1)$.

The proof of the case $p=+\infty$ is similar to that for $p\in[1,+\infty)$, so we omit it. 
\end{proof}

For illustration, the table below shows, for several specific values of $l$, the minimum number of norming constants required for uniqueness in the half-inverse problems considered in Theorem \ref{7.2} with $L^1(0,1)$ potentials.
\begin{center}
\begin{tabular}{c@{\qquad}cc}
\hline
$l$ &  $M^{(1)}_D(l)$ & $M^{(1)}_R(l)$\\
\hline
$-\frac12$ & $0$ & $0$\\
$0$        & $0$ & $0$\\
$1$        & $0$ & $1$\\
$2$        & $1$ & $1$\\
$3$        & $1$ & $2$\\
$4$        & $2$ & $2$\\
\hline
\end{tabular}
\end{center}

\vskip 1cm
\section{Proofs of Theorems \ref{8.4} and \ref{8.6}}\label{local}

This section is devoted to the proofs of Theorems \ref{8.4} and \ref{8.6}. Thus, we only consider $\beta=\infty$ and $l\in\mathbb{N}$. In this section one sets  
\begin{equation}\label{42}
M=M^{(2)}_D(l)=\left\lceil\frac{2l-3}{4}\right\rceil,
 \quad J\subset\mathbb{N},\quad \# J=M,
\end{equation}
where $M^{(2)}_D(l)$ is defined in \eqref{36}. Fix a model potential $q^*\in L^2(0,1)$. Recall the space $\mathcal X_{q^*}$ was defined in \eqref{41}, the maps $\mathcal S_J$ and $m(\cdot)$ were defined in \eqref{44} and \eqref{43}, respectively. 

The core step of this section is to prove that the Fréchet  derivative of $S_J$ at $q^*$ ( i.e., $D\mathcal S_J(q^*)$) is an isomorphism; see Lemma \ref{isom}. For this purpose, we will consider the map
\begin{equation}\label{dsq}
 D\mathcal S_J(q^*)[v] =\left(m(v),\{D\lambda_n(q^*)[v]-m(v)\}_{n\geq1},\{ D\tau_j(q^*)[v]\}_{j\in J}\right),\; v\in L^2(\tfrac{1}{2},1).
\end{equation}
Lemma \ref{8.1} will provide a certain relationship between $D\lambda_n(q^*)[v]$, $D\tau_j(q^*)[v]$ and $v$. Lemma \ref{8.3} will provide the Riesz basis property for a class of trigonometric functions.

To prove Theorems \ref{8.4} and \ref{8.6}, after adding the same sufficiently large constant to the model and comparison potentials, without loss of generality, we may assume that $\lambda_1(q^*)>0$. This leaves the norming constants and $d_J$ unchanged. Since $\lambda_1(q)$ is real analytic as a function of the $L^2$ potential, there exists a neighborhood of $q^*$ in $L^2(0,1)$ such that $\lambda_1(q)>0$ for all $q$ in this neighborhood.

We assume that $v\in L^2(\tfrac{1}{2},1)$. Without causing any ambiguity, we will also extend $v$ to the interval $[0,1]$ by setting $v=0$ on $[0,\tfrac{1}{2}]$, still denoted by $v$.
For $\varepsilon\in\mathbb{R}$, we denote 
$$q_\varepsilon:=q^*+\varepsilon v.$$
For simplicity, we will sometimes write $\lambda^*_n:=\lambda_n(q^*)$ and use the notation $\ip{\cdot}{\cdot}$ to mean the inner product on $L^2(0,\tfrac{1}{2})$. Recall that, for $n\in\mathbb{N}$ and $j\in J$, 
$$c_{\lambda^*_n}(y)=\cos(2\sqrt{\lambda^*_n}y),\quad d_{\lambda^*_j}(y) =\frac{y}{\sqrt{\lambda^*_j}}\sin(2\sqrt{\lambda^*_j}y).$$

\begin{lemma}\label{8.1}
There exists an isomorphism $\cA: L^2(\tfrac{1}{2},1)\rightarrow L^2(0,\tfrac{1}{2})$, such that for every $v\in L^2(\tfrac{1}{2},1)$, one has 

(1)
\begin{equation}\label{63}
 \ip{\cA v}{1}=m(v); 
\end{equation}

(2) 
\begin{equation}\label{66}
 \ip{\cA v}{ c_{\lambda^*_n}}=m(v)-a_n D\lambda_n(q^*)[v],
\end{equation}
where 
\begin{equation}\label{65}
a_n:=2\lambda_n^*\kappa_n^2(q^*)\tau_n(q^*),
\end{equation}
and for large $n$,  
\begin{equation}\label{67}
 a_n=1+O(n^{-1}); 
\end{equation}

(3) \begin{equation}\label{75}
 \ip{\cA v}{d_{\lambda^*_j}}=p_jD\tau_j(q^*)[v]+b_j D\lambda_j(q^*)[v],
\end{equation}
where
\begin{equation}\label{76}
p_j=-2\lambda^*_j\kappa_j^2(q^*)\neq0,
\end{equation}
and
\begin{equation}\label{77}
b_j=2\kappa_j^2(q^*)\tau_j(q^*)-2\lambda^*_j\kappa_j(q^*)\ddot\Delta(\lambda_j^*,q^*).
\end{equation}
\end{lemma}

\begin{proof}
The proof can be divided into three steps. 

\textbf{Step 1. The construction of $\cA$.} 
Since $\beta=\infty$, the solution defined in \eqref{8} satisfies $\psi(1,\lambda,q)=0$ and $\psi'(1,\lambda,q)=-1$. As in  \eqref{15}, one sets    
\begin{equation}\label{52}
 F_\varepsilon(\lambda):=W\left(\psi(\cdot,\lambda, q_\varepsilon), \psi(\cdot,\lambda, q^*)\right)\left(\tfrac{1}{2}\right).
\end{equation}
The Lagrange identity gives
\begin{equation}\label{53}
 F_\varepsilon(\lambda)=\varepsilon\int_{1/2}^1v(x) \psi(x,\lambda, q_\varepsilon) \psi(x,\lambda, q^*)\,\dd x.
\end{equation}
Let 
\begin{equation}\label{54}
 Q_{\varepsilon}(y):=\frac{l(l+1)}{(1-y)^2}+q_\varepsilon(1-y),\quad
 u_\varepsilon(y,\lambda):=\psi(1-y,\lambda, q_\varepsilon),\quad y\in \left[0,\tfrac{1}{2}\right].
\end{equation}
Then $u_\varepsilon(y,\lambda)$ is the solution to  
$$-u_\varepsilon''+Q_\varepsilon u_\varepsilon=\lambda u_\varepsilon,
 \quad
u_\varepsilon(0,\lambda)=0,
 \quad u_\varepsilon'(0,\lambda)=1.$$ 
Using \eqref{53} and the substitution $x=1-y$, we obtain
$$F_\varepsilon(\lambda)=\varepsilon\int_0^{1/2} (\mathcal{R}v)(y)u_\varepsilon(y,\lambda)u_0(y,\lambda)\,\dd y,$$
where $\mathcal{R}: L^2(\tfrac{1}{2},1)\rightarrow L^2(0,\tfrac{1}{2})$ is the reflection operator defined by 
$$(\mathcal{R}v)(y)=v(1-y),\quad y\in \left[0,\tfrac{1}{2}\right].$$

Note that 
$$Q_\varepsilon-Q_0=\varepsilon v(1-y)=\varepsilon (\mathcal{R}v)(y),\quad y\in \left[0,\tfrac{1}{2}\right].$$
Applying Proposition \ref{4.1} to the pair of potentials $Q_\varepsilon, Q_0$ with $\mu=0$ in \eqref{22}, one sees that there exists an operator $A_0(Q_\varepsilon, Q_0)$ which is an isomorphism of $L^2(0,\tfrac{1}{2})$, such that 
\begin{equation}\label{57}
2\lambda F_\varepsilon(\lambda)=\varepsilon\int_0^{1/2} (\mathcal{R}v)(y)\,\dd y-\varepsilon\int_0^{1/2}
 \left[A_0(Q_\varepsilon, Q_0)\mathcal{R}v\right](y)\cos(2\sqrt\lambda\,y)\,\dd y.
\end{equation}
Moreover, for sufficiently small $|\varepsilon|$, there is a constant $D$ such that
\begin{equation}\label{56}
 \|Q_\varepsilon\|_{L^2(0,1/2)}+\|Q_0\|_{L^2(0,1/2)}\leq D.
\end{equation}
The kernel estimate in \cite[Lemma 4.1]{HK2010} gives uniformly under \eqref{56}
\begin{equation}\label{59}
 \|A_{0}(Q_\varepsilon, Q_0)-A_{0}(Q_0, Q_0)\|
 \leq C_D\|Q_\varepsilon-Q_0\|_{L^1(0,1/2)},
\end{equation}
where $C_D$ is a positive constant depending on $D$.

Note that $\int_0^{1/2} (\mathcal{R}v)(y)\,\dd y=\int_{1/2}^1v(x)\,\dd x=m(v)$. We divide \eqref{57} by $\varepsilon$ and use \eqref{59}, then   letting $\varepsilon\rightarrow0$ yields 
\begin{equation}\label{62}
 \int_0^{1/2} (\cA v)(y)\cos(2\sqrt\lambda\,y)\,\dd y
 =m(v)-2\lambda G_v(\lambda),
\end{equation}
where 
\begin{equation}\label{61}
\cA:=A_0(Q_0, Q_0)\mathcal{R},\quad
G_v(\lambda):= \left.\partial_\varepsilon F_\varepsilon(\lambda)\right|_{\varepsilon=0}.
\end{equation}
Since $A_0(Q_0, Q_0)$ and $\mathcal{R}$ are isomorphisms, the operator $\cA$ is an isomorphism.

\textbf{Step 2. Proofs of (1) and (2).} 
Substituting $\lambda=0$ into \eqref{62}, we obtain \eqref{63}. Since $q_\varepsilon=q^*$ a.e. on $(0,\tfrac{1}{2})$, for the eigenvalue $\lambda_n^*$ we have 
\begin{equation}\label{b68}
F_\varepsilon(\lambda_n^*) =\kappa_n(q^*)\Delta(\lambda_n^*,q_\varepsilon).
\end{equation}
Since $\Delta(\lambda_n(q_\varepsilon), q_\varepsilon)\equiv0$, differentiating with respect to $\varepsilon$ and setting $\varepsilon=0$, together with \eqref{12}, yields 
\begin{equation}\label{bb68}
\partial_q\Delta(\lambda_n^*,q^*)[v] =-\dot\Delta(\lambda_n^*,q^*)D\lambda_n(q^*)[v]=\kappa_n(q^*)\tau_n(q^*) D\lambda_n(q^*)[v].
\end{equation}
Differentiating both sides of \eqref{b68} at $\varepsilon=0$, and using \eqref{61} and \eqref{bb68}, we obtain 
\begin{equation}\label{68}
 G_v(\lambda_n^*)=\kappa_n^2(q^*)\tau_n(q^*) D\lambda_n(q^*)[v].
\end{equation}
Putting $\lambda=\lambda^*_n$ in \eqref{62} yields \eqref{66}.

The estimate (2.17) in \cite{Serier2007} gives 
\begin{equation}\label{69}
\tau_n(q^*)=\frac{1}{2(\lambda_n^*)^{l+1}}\bigl(1+O(n^{-1})\bigr) \; \text{for large}\; n.
\end{equation}
On the other hand, the estimates (2.20) and (2.25) in \cite{KST2010}, together with \eqref{45}, give   
\begin{equation}\label{71}
\kappa^2_n(q^*)=\left|\phi'(1,\lambda_n^*,q^*)\right|^{-2}=(\lambda_n^*)^{l}\bigl(1+O(n^{-1})\bigr) \; \text{for large}\; n.
\end{equation}
One can obtain \eqref{67} from \eqref{69} and \eqref{71}.

\textbf{Step 3. Proof of (3).}
We fix $j\in J$ and use a dot to denote differentiation with respect to the spectral parameter $\lambda$. 
Differentiating \eqref{62} with respect to $\lambda$ and using \eqref{68}, we obtain 
\begin{equation}\label{bb77}
\begin{split}
\ip{\cA v}{d_{\lambda_j^*}}& =2G_v(\lambda_j^*)+2\lambda^*_j\dot G_v(\lambda_j^*) \\
   & =2\kappa_j^2(q^*)\tau_j(q^*) D\lambda_j(q^*)[v]+2\lambda^*_j\dot G_v(\lambda_j^*). 
\end{split}
\end{equation}

Differentiating \eqref{52} with respect to $\varepsilon$ at $\varepsilon=0$ and noting that $G_v(\lambda)= \left.\partial_\varepsilon F_\varepsilon(\lambda)\right|_{\varepsilon=0}$, one has 
\begin{equation}\label{ste1}
G_v(\lambda)=W\left(\partial_q\psi(\cdot,\lambda, q^*)[v], \psi(\cdot,\lambda, q^*)\right)\left(\tfrac{1}{2}\right).
\end{equation}
For each fixed $\lambda$, the equality $q_\varepsilon=q^*$ a.e. on $(0,\tfrac12)$ implies that $\phi(x,\lambda,q_\varepsilon)$ and $\phi'(x,\lambda,q_\varepsilon)$ are independent of $\varepsilon$ for $0<x\le\tfrac12$. Using \eqref{9} we obtain 
\begin{equation}\label{ste2}
\partial_q\Delta(\lambda, q^*)[v]=W\left(\partial_q\psi(\cdot,\lambda, q^*)[v], \phi(\cdot,\lambda, q^*)\right)\left(\tfrac{1}{2}\right).
\end{equation}
Differentiating both sides of the identity 
$$\psi(x,\lambda_j(q_\varepsilon), q_\varepsilon)=\kappa_j(q_\varepsilon)\phi(x,\lambda_j(q_\varepsilon), q_\varepsilon),\quad x\in(0,\tfrac{1}{2}],$$
with respect to $\varepsilon$ at $\varepsilon=0$, we have 
\begin{equation}\label{ste3}
\begin{split}
\partial_q\psi(x,\lambda_j^*, q^*)[v] =& \phi(x,\lambda_j^*, q^*)D\kappa_j(q^*)[v] \\
     & +\left(\kappa_j(q^*)\dot{\phi}(x,\lambda_j^*, q^*)-\dot{\psi}(x,\lambda_j^*, q^*)\right)D\lambda_j(q^*)[v].
\end{split}
\end{equation}

Differentiating the two Wronskian expressions \eqref{ste1} and \eqref{ste2} with respect to \(\lambda\), using $\psi(x,\lambda_j^*,q^*)=\kappa_j(q^*)\phi(x,\lambda_j^*,q^*)$ and \eqref{ste3}, we obtain
\begin{equation}\label{ste4}
\begin{split}
    & \dot{G}_v(\lambda_j^*)-\kappa_j(q^*)\partial_q\dot{\Delta}(\lambda_j^*, q^*)[v] \\
     & =W\left(\partial_q\psi(\cdot,\lambda_j^*, q^*)[v], \;\dot{\psi}(\cdot,\lambda_j^*, q^*)-\kappa_j(q^*)\dot{\phi}(\cdot,\lambda_j^*, q^*)\right)\left(\tfrac{1}{2}\right)\\
     & =D\kappa_j(q^*)[v]W\left(\phi(\cdot,\lambda_j^*, q^*), \;\dot{\psi}(\cdot,\lambda_j^*, q^*)-\kappa_j(q^*)\dot{\phi}(\cdot,\lambda_j^*, q^*)\right)\left(\tfrac{1}{2}\right).
\end{split}
\end{equation}
On the other hand, 
\begin{equation}\label{ste5}
\begin{split}
\dot{\Delta}(\lambda_j^*, q^*) & =W\left(\dot{\psi}(\cdot,\lambda_j^*, q^*), \phi(\cdot,\lambda_j^*, q^*)\right)+W\left(\psi(\cdot,\lambda_j^*, q^*), \dot{\phi}(\cdot,\lambda_j^*, q^*)\right)\\
 & =-W\left(\phi(\cdot,\lambda_j^*, q^*), \dot{\psi}(\cdot,\lambda_j^*, q^*)\right)+\kappa_j(q^*)W\left(\phi(\cdot,\lambda_j^*, q^*), \dot{\phi}(\cdot,\lambda_j^*, q^*)\right)\\
& =-W\left(\phi(\cdot,\lambda_j^*, q^*), \dot{\psi}(\cdot,\lambda_j^*, q^*)-\kappa_j(q^*)\dot{\phi}(\cdot,\lambda_j^*, q^*)\right).   
\end{split}
\end{equation}
Then \eqref{ste4} and \eqref{ste5} yield
\begin{equation}\label{ste6}
\dot{G}_v(\lambda_j^*)=\kappa_j(q^*)\partial_q\dot{\Delta}(\lambda_j^*, q^*)[v]-\dot{\Delta}(\lambda_j^*, q^*)D\kappa_j(q^*)[v]. 
\end{equation}

Finally, the identity \eqref{12} gives 
\begin{equation}\label{ste7}
\dot\Delta(\lambda_j( q_\varepsilon),  q_\varepsilon) =-\kappa_j( q_\varepsilon)\tau_j( q_\varepsilon).
\end{equation}
Differentiating \eqref{ste7} with respect to $\varepsilon$ at $\varepsilon=0$ yields 
\begin{equation}\label{ste8}
 \partial_q\dot{\Delta}(\lambda_j^*,q^*)[v]+\ddot{\Delta}(\lambda_j^*,q^*)D\lambda_j(q^*)[v] =-\tau_j(q^*)D\kappa_j(q^*)[v]-\kappa_j(q^*)D\tau_j(q^*)[v].
\end{equation}
From \eqref{ste6} and \eqref{ste8}, and using \eqref{12} again, we obtain
\begin{equation}\label{8.50}
 \dot G_v(\lambda_j^*) =-\kappa_j^2(q^*)D\tau_j(q^*)[v] -\kappa_j(q^*)\ddot\Delta(\lambda_j^*,q^*)D\lambda_j(q^*)[v]. 
\end{equation} 
Then \eqref{bb77} together with \eqref{8.50} implies \eqref{75}. 
\end{proof}

\begin{lemma}\label{8.3}
The family 
\begin{equation}\label{83}
\mathcal E_*:=\{1\}\cup\{ c_{\lambda^*_n}| n\geq1\}\cup\{d_{\lambda^*_j}|j\in J\}
\end{equation}
is a Riesz basis in $L^2(0,\tfrac{1}{2})$.
\end{lemma}

\begin{proof}
We first prove the completeness. Theorem \ref{7.2} with $p=2$ and the selected index set $J$ yields the uniqueness property since $\# J=M_D^{(2)}(l)=M$.  Since $0\notin\{\lambda_n^*|n\geq1\}$, we may choose $\mu=0$. Therefore the necessary and sufficient condition in Theorem \ref{5.4} implies that $\mathcal E_*$ is closed, equivalently, complete for $p=2$, in $L^2(0,\tfrac{1}{2})$. 

On the other hand, write $l=2M$ or $l=2M+1$ and set
\begin{equation}\label{85}
 e_n(y):=
 \begin{cases}
 \cos(2(n+M)\pi y),&l=2M,\\
 \cos((2(n+M)+1)\pi y),&l=2M+1.
 \end{cases}
\end{equation}
From \eqref{45} one has 
\begin{equation}\label{84}
\sqrt{\lambda^*_n}=\pi\left(n+\frac l2\right)+O(n^{-1}).
\end{equation}
It is easy to conclude by using \eqref{84} that 
\begin{equation}\label{86}
 \sum_{n=1}^{\infty} \| c_{\lambda^*_n}-e_n\|_{L^2(0,1/2)}^2<\infty.
\end{equation}
Set
\begin{equation*}
E_l:=
\begin{cases}
\left\{1,\cos(2\pi y),\cdots,\cos(2M\pi y)\right\}, &   l=2M, \\
\left\{\cos(\pi y),\cos(3\pi y),\cdots,\cos((2M+1)\pi y)\right\}, &  l=2M+1.
\end{cases}
\end{equation*}
Then $E_l\cup\{e_n\}_{n\ge1}$ is an orthogonal basis of $L^2(0,\frac12)$. Since $\{1\}\cup\{d_{\lambda^*_j}:j\in J\}$ and $E_l$ are two finite sets with the same number of elements, the inequality \eqref{86} implies that $\mathcal E_*$ is \emph{quadratically close} (\cite[Assumption 1.8.1]{FY2001}) to $E_l\cup\{e_n\}_{n\ge1}$. From Proposition 1.8.5 in \cite{FY2001}, one sees that $\mathcal E_*$ is a Riesz basis. 
\end{proof}

\begin{lemma}\label{isom}
The map 
\begin{equation*}
 D\mathcal S_J(q^*):L^2(\tfrac{1}{2},1)\longrightarrow\mathcal Y
\end{equation*}
is an isomorphism.
\end{lemma}

\begin{proof}
Define 
$$\mathcal{U}f:=\left(\ip{f}{1},\left\{\ip{f}{ c_{\lambda^*_n}}\right\}_{n\in\N},\left\{\ip{f}{d_{\lambda^*_j}}\right\}_{j\in J}\right),\quad f\in L^2(0,\tfrac{1}{2}). $$
Lemma \ref{8.3} implies that $\mathcal{U}:L^2(0,\tfrac{1}{2})\longrightarrow\mathcal Y$ is an isomorphism.

For $n\in\mathbb{N}$ and $j\in J$, denote 
\begin{equation}\label{72}
m_0:=m(v),\quad \xi_n:= D\lambda_n(q^*)[v]-m(v), \quad \eta_j:=D\tau_j(q^*)[v].
\end{equation}
From \eqref{dsq} and \eqref{72} one has 
\begin{equation}\label{90}
 D\mathcal S_J(q^*)[v] =\bigl(m_0,\{\xi_n\}_{n\geq1},\{\eta_j\}_{j\in J}\bigr).
\end{equation}
Now \eqref{66} and \eqref{75} become, respectively, 
\begin{equation}\label{73}
 \ip{\cA v}{ c_{\lambda^*_n}}=-a_n \xi_n+(1-a_n)m_0,
\end{equation}
and
\begin{equation}\label{bb73}
  \ip{\cA v}{d_{\lambda^*_j}}=p_j \eta_j+b_j (\xi_j+m_0).
\end{equation}

Define the triangular operator $\mathfrak{T}:\mathcal Y\to\mathcal Y$ by
\begin{equation}\label{91}
\mathfrak{T}(m_0,\xi,\eta)=\left( m_0,\{-a_n \xi_n+(1-a_n)m_0\}_{n\in\N},\{p_j \eta_j+b_j(\xi_j+m_0)\}_{j\in J}\right),
\end{equation}
with $m_0\in\mathbb{R}$, $\xi=\{\xi_n\}_{n\in\N}\in \ell^2(\N)$ and $\eta=\{\eta_j\}_{j\in J}\in\mathbb{R}^M$.
Equations \eqref{63}, \eqref{73}, and \eqref{bb73} now give the operator identity
\begin{equation}\label{94}
\mathcal{U}\cA=\mathfrak{T}\,D\mathcal S_J(q^*).
\end{equation}
By \eqref{65}, $a_n>0$ for every $n$, then the estimates \eqref{67} and the fact $p_j\neq0$ show that $\mathfrak{T}$ is boundedly invertible.  Indeed, suppose that  
$$\mathfrak{T}(m_0,\xi,\eta)=(m_0,\tilde{\xi}, \tilde{\eta})$$
with $\tilde{\xi}=\{\tilde\xi_n\}_{n\in\N}\in \ell^2(\N)$ and $\tilde{\eta}=\{\tilde\eta_j\}_{j\in J}\in\mathbb{R}^M$. Then the inverse of $\mathfrak{T}$ is given by
\begin{align}
 \xi_n&=-a_n^{-1}\bigl(\tilde{\xi}_n-(1-a_n)m_0\bigr),
 \label{92}\\
 \eta_j&=p_j^{-1}\bigl(\tilde{\eta}_j-b_j(\xi_j+m_0)\bigr).
 \label{93}
\end{align}
Furthermore, it is also easy to verify from \eqref{92} and \eqref{93} that $\mathfrak{T}^{-1}$ is bounded. 
Since $\cA$, $\mathcal{U}$, and $\mathfrak{T}$ are all isomorphisms, the operator identity \eqref{94} implies that $D\mathcal S_J(q^*)$ is an isomorphism.   
\end{proof}

Using Lemma \ref{isom}, now we prove Theorem \ref{8.4}.

\begin{proof}[Proof of Theorem \ref{8.4}]
Using Lemma \ref{isom} we may define  
\begin{equation}\label{96}
 \gamma:=\|D\mathcal S_J(q^*)^{-1}\|^{-1}>0,
\end{equation}
and thus 
\begin{equation}\label{97}
 \|D\mathcal S_J(q^*)v\|_{\mathcal Y}\geq\gamma\|v\|_{L^2(1/2,1)},\quad v\in L^2(\tfrac{1}{2},1).
\end{equation}
Since the map $\mathcal S_J$ is real analytic, the derivative $D\mathcal S_J(q)$ is continuous in the operator norm. We can choose $C_*>0$ such that
\begin{equation}\label{100}
 \sup_{\|q-q^*\|_{L^2(1/2,1)}<C_*}
 \|D\mathcal S_J(q)-D\mathcal S_J(q^*)\|
 \leq\frac\gamma2.
\end{equation}
For $v=q_2-q_1$ and $q_t=q_1+tv$, we have 
\begin{equation}\label{bb100}
\mathcal S_J(q_2)-\mathcal S_J(q_1)=\int_{0}^{1}D\mathcal S_J(q_t)v\,\dd t.
\end{equation}
Substituting the identity 
$$D\mathcal S_J(q_t)v=D\mathcal S_J(q^*)v+\left(D\mathcal S_J(q_t)-D\mathcal S_J(q^*)\right)(v)$$
into the integral term of \eqref{bb100}, we obtain
\begin{equation}\label{bbb100}
\mathcal S_J(q_2)-\mathcal S_J(q_1)=D\mathcal S_J(q^*)v+\int_0^1\bigl(D\mathcal S_J(q_t)-D\mathcal S_J(q^*)\bigr)v\,\dd t.
\end{equation}
By considering the construction of $ d_J$ in \eqref{49}, and using \eqref{97}, \eqref{100}, and \eqref{bbb100}, we obtain 
\begin{equation}\label{b100}
d_J\left(\mathcal S_J(q_1),\mathcal S_J(q_2)\right)=\|\mathcal S_J(q_2)-\mathcal S_J(q_1)\|_{\mathcal Y}\geq\frac\gamma2\|q_2-q_1\|_{L^2(1/2,1)}.
\end{equation}
Moreover, the estimate \eqref{100} gives the explicit uniform upper bound
\begin{equation}\label{bbbb100}
\sup_{\|q-q^*\|_{L^2(1/2,1)}<C_*}\|D\mathcal S_J(q)\| \leq \|D\mathcal S_J(q^*)\|+\frac\gamma2.
\end{equation}
From \eqref{bb100} and \eqref{bbbb100}, we obtain  
\begin{equation}\label{bbbbb100}
\begin{aligned}
d_J\left(\mathcal S_J(q_1),\mathcal S_J(q_2)\right) =&\|\mathcal S_J(q_2)-\mathcal S_J(q_1)\|_{\mathcal Y}\\
 &\leq\int_0^1\|D\mathcal S_J(q_t)\|
              \|q_2-q_1\|_{L^2(1/2,1)}\,\dd t\\
 &\leq \left(\|D\mathcal S_J(q^*)\|+\frac\gamma2\right) \|q_2-q_1\|_{L^2(1/2,1)}.
 \end{aligned}
\end{equation}
Taking \eqref{b100} and \eqref{bbbbb100} into account, we can take $C_1=\gamma/2$ and $C_2=\|D\mathcal S_J(q^*)\|+\gamma/2$. The proof is complete. 
\end{proof}

\begin{proof}[Proof of Theorem \ref{8.6}]
Let $C_*$ and $C_1$ be the constants in Theorem \ref{8.4}. By Lemma \ref{isom}, $D\mathcal S_J(q^*)$ is a bounded isomorphism onto $\mathcal Y$. Since $\mathcal S_J$ is real analytic and $\mathcal S_J(q^*)=0$, the real-analytic inverse function theorem yields an open neighborhood $U_1$ of $q^*$ in $\mathcal X_{q^*}$ and an open neighborhood $U_2$ of $0$ in $\mathcal Y$ such that 
$$\mathcal S_J: U_1\longrightarrow U_2$$
has a real-analytic inverse. Shrinking these neighborhoods if necessary, we may assume that
$$U_1\subset\left\{q\in\mathcal X_{q^*}|\|q-q^*\|_{L^2(1/2,1)}<C_*\right\}.$$
Choose $\delta>0$ such that 
$$B_{\mathcal Y}(0,\delta):=\{\mathfrak{s}\in\mathcal Y|\|\mathfrak{s}\|_{\mathcal Y}<\delta\}\subset U_2.$$
For any $\mathfrak{s}\in B_{\mathcal Y}(0,\delta)$, define
$$ q:=\mathcal S_J^{-1}(\mathfrak{s}).$$
Then $\mathcal S_J(q)=\mathfrak{s}$. 

Uniqueness follows from Theorem 1.2 with $p=2$, since any two potentials in $\mathcal X_{q^*}$ with the same image under $\mathcal S_J$ share the full spectrum and the selected norming constants.

Finally, applying Theorem \ref{8.4} to $q$ and $q^*$, we obtain
\begin{equation*}
C_1\|q-q^*\|_{L^2(1/2,1)}\leq \|\mathcal S_J(q)-\mathcal S_J(q^*)\|_{\mathcal Y}=\|\mathfrak s\|_{\mathcal Y}.
\end{equation*}
Thus the required estimate \eqref{th162} holds with $C=C_1^{-1}$.
\end{proof}

\begin{remark}\label{8.5}
(1) Theorem \ref{8.6} also yields the conclusion of Theorem \ref{7.3} for positive integer $l$ and $p=2$. (2) Extensions of Theorems \ref{8.4} and \ref{8.6} to noninteger $l$ or to the Robin boundary condition $\beta\in\R$ require additional nonharmonic frame analysis and are not addressed here.  
\end{remark}

\vskip 1cm
\appendix
\section{Two auxiliary lemmas}\label{appe}

In this section we provide two auxiliary lemmas. First, we present a slightly modified version of \cite[Lemma 2.1]{Horvath2005}, in order to avoid the extreme case where $\|q_0\|=0$ in equation (2.2) of Lemma 2.1 in \cite{Horvath2005}. 
\begin{lemma}[\cite{Horvath2005}]\label{h2005}
Let $B_1$ and $B_2$ be Banach spaces. For every $Q\in B_1$, let
$$A_Q:B_1\longrightarrow B_2$$
be a bounded linear operator. Suppose that, for some $Q_0\in B_1$, the operator $A_{Q_0}$ is an isomorphism. Assume further that there exist constants $\rho>0$ and $C>0$ such that
$$\|A_{Q_1}-A_{Q_2}\|\leq C\|Q_1-Q_2\|_{B_1}$$
for all $Q_1,Q_2\in B_1$ satisfying
$$\|Q_i-Q_0\|_{B_1}\leq\rho, \quad i=1,2.$$
Then the set
$$\bigl\{A_Q(Q-Q_0):Q\in B_1\bigr\}$$
contains an open ball in $B_2$ centered at the origin.
\end{lemma}

Fix $b>0$, and write $\|\cdot\|_p$ for the norm in $L^p(0,b)$.  For $Q\in L^p(0,b)$ and $\beta\in\mathbb R$, let $u_{Q,\beta}(y,z)$ solve
\begin{equation}\label{aa1}
\begin{gathered}
-u_{Q,\beta}''(y,z)+Q(y)u_{Q,\beta}(y,z)=z^2u_{Q,\beta}(y,z),\\
u_{Q,\beta}(0,z)=1,\quad u_{Q,\beta}'(0,z)=\beta.
\end{gathered}
\end{equation}
The following lemma is used in the proofs of Theorem \ref{7.0} and Proposition \ref{4.1}.

\begin{lemma}\label{a.1}
Fix $Q_0\in L^p(0,b)$, $p\in [1,+\infty]$, $\beta_0\in\mathbb R$, and suppose that $Q, \tilde{Q}\in L^p(0,b)$, $\beta, \tilde{\beta}\in\mathbb R$.

(1) There is a continuous kernel $P_{Q,\beta,\tilde{Q},\tilde{\beta}}$ on $\Omega:=\{(y,t)| 0\leq t\leq y\leq b\}$, 
independent of $z$, such that, for $z\in\mathbb C$,
\begin{equation}\label{aa2}
u_{Q,\beta}(y,z)u_{\tilde{Q},\tilde{\beta}}(y,z)=\frac12+\frac12\cos(2zy)+\int_0^y P_{Q,\beta,\tilde{Q},\tilde{\beta}}(y,t)\cos(2zt)\,\dd t;
\end{equation}

(2) Define the bounded operator $B(Q,\beta,\tilde{Q},\tilde{\beta}): L^p(0,b)\rightarrow L^p(0,b)$ by
\begin{equation}\label{aa3}
 [B(Q,\beta,\tilde{Q},\tilde{\beta})f](t) :=\frac12 f(t)+\int_t^b P_{Q,\beta,\tilde{Q},\tilde{\beta}}(y,t)f(y)\,\dd y.
\end{equation}
Then for every $D>0$, there is a constant $C=C(\|Q_0\|_p,|\beta_0|,D,b)$ such that
\begin{equation}\label{aa4}
 \|B(Q_0,\beta_0,Q,\beta)-B(Q_0,\beta_0,\tilde{Q},\tilde{\beta})\|\leq C\bigl(\|Q-\tilde Q\|_p+|\beta-\tilde\beta|\bigr),
\end{equation}
whenever 
\begin{equation}\label{bound1}
\|Q-Q_0\|_p+|\beta-\beta_0|\leq D,\quad \|\tilde Q-Q_0\|_p+|\tilde\beta-\beta_0|\leq D.
\end{equation}
\end{lemma}

\begin{proof}

For $Q\in L^p(0,b)$, we consider the solution $e_Q(y,z)$ of 
$$-e_Q''+Qe_Q=z^2e_Q,\quad e_Q(0,z)=1,\quad e_Q'(0,z)={\rm i}z,\quad {\rm i}=\sqrt{-1}.$$
Extending $Q$ by zero to $(-b,0)$ and applying Lemma 5.1 in \cite{Horvath2005}, we obtain the following representation for this solution:
\begin{equation}\label{aa5}
 e_Q(y,z)=\mathrm e^{\mathrm izy}+\int_{-y}^y L_Q(y,t)\mathrm e^{\mathrm izt}\,\dd t,\quad 0\leq y\leq b,  
\end{equation}
where $L_Q(y,t)$ is continuous on the region $\Omega^*:=\{(y,t)| 0\leq |t|\leq y\leq b\}$.  For every fixed $E>0$, if $\|Q\|_p,\|\tilde Q\|_p\leq E$, the kernel estimates on the  region $\Omega^*$ are
\begin{equation}\label{aa6}
\|L_Q\|_{L^{\infty}(\Omega^*)}\leq C_E,\quad \|L_Q-L_{\tilde Q}\|_{L^{\infty}(\Omega^*)} \leq C_E\|Q-\tilde Q\|_p,
\end{equation}
with $C_E$ depending only on $E$ and $b$.

Define 
$$c_Q(y,z):=\frac{e_Q(y,z)+e_Q(y,-z)}2,\quad s_Q(y,z):=\frac{e_Q(y,z)-e_Q(y,-z)}{2\mathrm iz}.$$
The value of $s_Q(y, z)$ at $z=0$ is defined by continuity. Then $(c_Q(0,z), c'_Q(0,z))=(1,0)$, $(s_Q(0,z), s'_Q(0,z))=(0,1)$, and $u_{Q,\beta}=c_Q+\beta s_Q$. 
From \eqref{aa5} and using
$$\frac{\sin(zs)}{z}=\int_0^s\cos(zt)\,\dd t,$$
by Fubini's theorem one has 
\begin{equation}\label{aa7}
u_{Q,\beta}(y,z)=\cos(zy)+\int_0^y K_{Q,\beta}(y,t)\cos(zt)\,\dd t,
\end{equation}
where
\begin{equation}\label{aa8}
K_{Q,\beta}(y,t)=L_Q(y,t)+L_Q(y,-t)+\beta\left[1+\int_t^y\bigl(L_Q(y,s)-L_Q(y,-s)\bigr)\,\dd s\right].
\end{equation}

For clarity, put
\begin{equation}\label{lj}
L_Q^+(y,t):=L_Q(y,t)+L_Q(y,-t),\quad J_Q(y,t):=1+\int_t^y \bigl(L_Q(y,s)-L_Q(y,-s)\bigr)\,\dd s.
\end{equation}
Thus 
\begin{equation}\label{lj1}
K_{Q,\beta}(y,t)=L_Q^+(y,t)+\beta J_Q(y,t),
\end{equation}
and
\begin{equation}\label{lj2}
\begin{split}
K_{Q,\beta}(y,t)-K_{\tilde Q,\tilde \beta}(y,t)=& L_Q^+(y,t)-L_{\tilde Q}^+(y,t)+(\beta-\tilde \beta)J_Q(y,t) \\
     & +\tilde\beta\left(J_Q(y,t)-J_{\tilde Q}(y,t)\right).
\end{split}
\end{equation}
The bounds in \eqref{aa6} and the definition of $J_Q$ in \eqref{lj} imply
\begin{equation}\label{lj3}
\|J_Q\|_{L^{\infty}(\Omega)}\leq1+2bC_E,\quad \|J_Q-J_{\tilde Q}\|_{L^{\infty}(\Omega)}\leq2bC_E\|Q-\tilde Q\|_p.
\end{equation}
Consequently, on every set of $Q$, $\tilde Q$, $\beta$, $\tilde \beta$, such that $\|Q\|_p,\|\tilde Q\|_p,|\beta|,|\tilde \beta|$ are bounded, the kernels $K_{Q,\beta}$ are uniformly bounded. Moreover, from \eqref{aa6}, \eqref{lj2}, and \eqref{lj3},  one sees that there exists $C_3>0$, which depends on  $\|Q\|_p,\|\tilde Q\|_p,|\beta|,|\tilde \beta|,$ and $b$, such that 
\begin{equation}\label{aa9}
\|K_{Q,\beta}-K_{\tilde Q,\tilde \beta}\|_{L^{\infty}(\Omega)}\leq C_3\bigl(\|Q-\tilde Q\|_p+|\beta-\tilde \beta|\bigr).
\end{equation}

Multiplying the representations \eqref{aa7} for $u_{Q,\beta}$ and $u_{\tilde Q,\tilde\beta}$, and changing variables in the integrals yield \eqref{aa2}, with  
\begin{equation}\label{aa11}
 \begin{split}
P_{Q,\beta,\tilde{Q},\tilde{\beta}}(y,t) ={}&K_{Q,\beta}(y,|2t-y|)+K_{\tilde{Q},\tilde{\beta}}(y,|2t-y|)\\
&+\int_{\max(0,2t-y)}^{\min(y,2t)}K_{Q,\beta}(y,s)K_{\tilde{Q},\tilde{\beta}}(y,2t-s)\,\dd s+R(y,t),
 \end{split}
\end{equation}
where
$$R(y,t)=\int_0^{y-2t}\!\bigl[K_{Q,\beta}(y,s+2t)K_{\tilde{Q},\tilde{\beta}}(y,s)+K_{Q,\beta}(y,s)K_{\tilde{Q},\tilde{\beta}}(y,s+2t)\bigr]\, \dd s$$
for $0\leq2t\leq y$, and $R(y,t)=0$ for $y<2t\leq2y$. From \eqref{aa8} and \eqref{aa11}, one can see that $P_{Q,\beta,\tilde{Q},\tilde{\beta}}$ is continuous on $\Omega$. This proves (1).

The continuity of $P_{Q,\beta,\tilde{Q},\tilde{\beta}}$ and Hölder's inequality imply that $B(Q,\beta,\tilde{Q},\tilde{\beta})$ is bounded on $L^p(0,b)$. Moreover, 
$$\|B(Q,\beta,\tilde Q,\tilde\beta)\|\leq\frac12+b\|P_{Q,\beta,\tilde Q,\tilde\beta}\|_{L^{\infty}(\Omega)}.$$
The inequalities in \eqref{bound1} imply 
$${\rm max}\left\{\|Q\|_p,\ \|\tilde Q\|_p\right\}\leq \|Q_0\|_p+D,\quad {\rm max}\left\{|\beta|,\ |\tilde\beta|\right\}\leq |\beta_0|+D.$$
Thus the preceding kernel estimates hold uniformly on this set, with constants depending only on $\|Q_0\|_p$, $|\beta_0|$, $D$, and $b$. Now, the relations \eqref{aa9} and \eqref{aa11} give 
\begin{equation}\label{aa12}
 \begin{split}
 \|P_{Q_0,\beta_0,Q,\beta}-P_{Q_0,\beta_0,\tilde{Q},\tilde{\beta}}\|_{L^{\infty}(\Omega)}
 &\leq\bigl(1+3b\|K_{Q_0,\beta_0}\|_{L^{\infty}(\Omega)}\bigr)
       \|K_{Q,\beta}-K_{\tilde Q,\tilde\beta}\|_{L^{\infty}(\Omega)}\\
 &\leq C_4\bigl(\|Q-\tilde Q\|_p+|\beta-\tilde\beta|\bigr)
 \end{split}
\end{equation}
for some $C_4>0$. Finally, for $f\in L^p(0,b)$, the definition in \eqref{aa3} and the inequality \eqref{aa12} give 
\begin{equation}\label{aa14}
\begin{aligned}
 &\|[B(Q_0,\beta_0,Q,\beta)-B(Q_0,\beta_0,\tilde{Q},\tilde{\beta})]f\|_p\\
 &\leq b\|P_{Q_0,\beta_0,Q,\beta}
-P_{Q_0,\beta_0,\tilde Q,\tilde\beta}\|_{L^{\infty}(\Omega)}
\|f\|_p\\
 &\leq C\bigl(\|Q-\tilde Q\|_p+|\beta-\tilde\beta|\bigr)\|f\|_p
 \end{aligned}
\end{equation}
for some $C=C(\|Q_0\|_p,|\beta_0|,D,b)>0$. The estimate \eqref{aa14} implies \eqref{aa4}. This proves~(2). 
\end{proof}

\bigskip

\noindent {\bf Acknowledgments.}
The authors gratefully acknowledge Professor Sergey Buterin for some valuable advice. 
The work of the author N. P. Bondarenko was supported by the Moscow Center for Fundamental and Applied Mathematics (grant of the Ministry of Education and Science of Russia No. 075-15-2025-345).
The work of the authors X. J. Xu and C. F. Yang was supported by the National Natural Science Foundation of China (No. 11871031 and No. 12501216) and by the Natural Science Foundation of the Jiangsu Province of China (BK 20201303 and BK 20241437). 

\vskip 0.5cm
\noindent {\bf Data availability.}
No data were used for the research described in the article.

\vskip 0.5cm
\noindent {\bf Conflict of Interest Statement.} 
The authors report no conflict of interest.

\end{document}